\documentstyle[11pt]{article}
       \font\tenmsb=msbm10
       \font\sevenmsb=msbm7
       \font\fivemsb=msbm5
       \catcode`\@=11
       \ifx\amstexloaded@\relax\catcode`\@=\active
       \endinput\else\let\amstexloaded@\relax\fi
       \def\spaces@{\space\space\space\space\space}
       \def\spaces@@{\spaces@\spaces@\spaces@\spaces@\spaces@}
       \def\space@.{\futurelet\space@\relax}
       \space@. %
       \def\Err@#1{\errhelp\defaulthelp@\errmessage{AmS-TeX error: #1}}
       \def\relaxnext@{\let\next\relax}
       \def\accentfam@{7}
       
       \def\noaccents@{\def\accentfam@{0}}
       \def\Cal{\relaxnext@\ifmmode\let\next\Cal@\else
       \def\next{\Err@{Use \string\Cal\space only in math mode}}\fi\next}
       \def\Cal@#1{{\Cal@@{#1}}}
       \def\Cal@@#1{\noaccents@\fam\tw@#1}

       \def\Bbb{\relaxnext@\ifmmode\let\next\Bbb@\else
       \def\next{\Err@{Use \string\Bbb\space only in math mode}}\fi\next}
       
       \def\Bbb@#1{{\Bbb@@{#1}}}
       \def\Bbb@@#1{\noaccents@\fam\msbfam#1}
       \newfam\msbfam
       \textfont\msbfam=\tenmsb
       \scriptfont\msbfam=\sevenmsb
       \scriptscriptfont\msbfam=\fivemsb

\def\N{{\Bbb N}}
\def\Z{{\Bbb Z}}
\def\R{{\Bbb R}}
\def\T{{\Bbb T}}
\def\C{{\Bbb C}}

\def\M{{\Bbb M}}

\newtheorem{Theorem}{Theorem}
\newtheorem{Lemma}{Lemma}[section]
\newtheorem{Proposition}{Proposition}
\newtheorem{Corollary}{Corollary}
\newtheorem{Remark}{Remark}[section]

\newtheorem{Definition}{Definition}[section]

\@addtoreset{equation}{section}

\newcommand{\qed}{\nolinebreak\hfill\rule{2mm}{2mm}
\par\medbreak}
\newcommand{\proof}{\par\medbreak\it Proof: \rm}

\newcommand{\la}{\langle }
\newcommand{\ra}{\rangle }

\newcommand{\beq}{\begin{equation} }
\newcommand{\eeq}{\end{equation} }

\begin{document}
\setlength{\columnsep}{5pt}
\title{Ballistic Transport in One-Dimensional Quasi-Periodic Continuous Schr\"odinger Equation}

\author{Zhiyan Zhao\footnote{This work is supported by ANR grant
``ANR-14-CE34-0002-01'' for the project ``Dynamics and CR geometry".}
\\
 {\footnotesize Laboratoire J.A. Dieudonn\'{e}, Universit\'{e} de Nice-Sophia Antipolis(Parc Valrose)}\\
{\footnotesize 06108 NICE Cedex 02, FRANCE}\\
 {\footnotesize Email: zhiyan.zhao@unice.fr}}

%\date{put in custom date; delete this line for today's date}
%\maketitle

\date{}
\maketitle

\begin{abstract}
For the solution $q(t)$ to the one-dimensional continuous Schr\"odinger equation
$${\rm i}\partial_t{q}(x,t)=-\partial_x^2 q(x,t) + V(\omega x) q(x,t), \quad x\in\R,$$
with $\omega\in\R^d$ satisfying a Diophantine condition, and $V$ a real-analytic function on $\T^d$, we consider the growth rate of the diffusion norm $\|q(t)\|_{D}:=\left(\int_{\R}x^2|q(x,t)|^2dx\right)^{\frac12}$ for any non-zero initial condition $q(0)\in H^{1}(\R)$ with $\|q(0)\|_D<\infty$.
We prove that $\|q(t)\|_{D}$ grows {\it linearly} with $t$ if $V$ is sufficiently small.
\end{abstract}

%\tableofcontents

\section{Introduction and main result}
\noindent

We consider the quasi-periodic Schr\"{o}dinger equation in one space dimension:
\begin{equation}\label{QP_Sch_EQ}
{\rm i}\partial_t{q}(x,t)=-\partial_x^2 q(x,t) + V(\omega x) q(x,t), \quad x\in\R,
\end{equation}
with $q(x,t)\in \C$ and $(x,t)\in\R\times \R$, where $V:\T^d\rightarrow \R$ is analytic in a complex neighbourhood of $\T^d$ $\{z\in\C^d: |\Im z|<r\leq 1\}$, and $\omega\in\R^d$ satisfies the Diophantine condition, i.e., there exist $\gamma>0$, $\tau>d-1$, such that
$$
\inf_{j\in\Z}\left|\frac{\la k, \omega\ra}{2}-j\pi\right| > \frac{\gamma}{|k|^{\tau}},\quad \forall \,  k\in\Z^d\setminus\{0\}.
$$
We would like to observe the growth rate with $t$ of the diffusion norm
$$\|q(t)\|_{D}:=\left(\int_{\R} x^2|q(x,t)|^2 dx\right)^{\frac12}$$
provided that $q(0)\in H^{1}(\R)\setminus\{0\}$ with $\|q(0)\|_D< \infty$.

The diffusion norm $\|\cdot\|_D$ is a weighted $L^2-$norm. It is well known that the $L^2-$norm $(\int_{\R}|q(x,t)|^2 dx)^\frac12$ is conserved for Eq.(\ref{QP_Sch_EQ}) since $e^{-{\rm i}tH}$ is unitary.
The initial condition $\|q(0)\|_D< \infty$ indicates the concentration on some range at the initial moment, and the diffusion norm $\|q(t)\|_D$ measures the propagation into the range where $|x|\gg 1$.

With the well-localized initial condition $\|q(0)\|_D< \infty$, we have that $\|q(t)\|_{D}<\infty$ for any finite $t$.
More precisely, for the bounded potential $V$ as in (\ref{QP_Sch_EQ}), we have the general ballistic upper bound, i.e., there exists a numerical constant $c>0$, such that
\begin{equation}\label{upper_bound}
\|q(t)\|_D\leq \|q(0)\|_D + c(\|q(0)\|_{H^1(\R)}+\|q(0)\|_{D}) t.
\end{equation}
(See Theorem 2.1 of \cite{RS}).
This is also related to the Lieb-Robinson bound\cite{LR}.

Since we are considering the linear equation (\ref{QP_Sch_EQ}), the behaviour of its solution is determined by the spectral property of the linear Schr\"odinger operator
\begin{eqnarray*}
H: \; {\cal C}_c^{\infty}(\R)  &\rightarrow& L^2(\R) \\
   q(x)  &\mapsto&  -q''(x)+V(\omega x) q(x) .
  \end{eqnarray*}
It has been shown in \cite{E92} that the spectrum of $H$ is purely absolutely continuous if $V$ is sufficiently small.
Inspired by RAGE Theorem\cite{CFKS}, it is natural to get the propagation which is related to the growth of $\|q(t)\|_D$ in this case.

More rigorously, one could expect ``ballistic motion" for Eq.(\ref{QP_Sch_EQ}) if the spectrum of the corresponding linear operator has
the absolutely continuous component.
Normally, it is interpreted as the linear growth with time of the diffusion norm.
A time-averaged statement by Guarneri-Combes-Last theorem\cite{Last} shows that, in the presence of absolutely continuous spectrum, we have
$$
\liminf_{t\to\infty}\frac1T\int_0^T \|q(t)\|_D \, dt\geq C T
$$
for some positive constant $C$.
A recent work by Karpeshina-Lee-Shterenberg-Stolz\cite{KLSS} shows the existence of ballistic transport for the Schr\"odinger equation with limit-periodic or quasi-periodic potential in two space dimension, under certain regularity assumptions on the potential which have been used
in prior work to establish the existence of an absolutely continuous component and other spectral properties.
Related to the lower bound of the Cesaro means of the diffusion norm for the solution to the continuous Schr\"odinger equation, we can also refer to \cite{AschKnauf, DamanikLenzStolz}.

\

In this paper, we try to go beyond the time-averaged version of ballistic transport in Guarneri-Combes-Last theorem.
For the solution to Eq.(\ref{QP_Sch_EQ}), we are going to show the linear growth of the diffusion norm.

\begin{Theorem}\label{thm-qp}
Consider the solution $q(t)$ of Eq.(\ref{QP_Sch_EQ}) with the non-zero initial condition $q(0)\in H^1(\R)$ and $\|q(0)\|_D<\infty$.
There exists an $\varepsilon_*=\varepsilon_*(\gamma,\tau,r)$ such that
if $|V|_r=\varepsilon_0<\varepsilon_*$, then there is a constant $0<C<\infty$, depending on $\varepsilon_0$ and $q(0)$, such that
$$ \liminf_{t\rightarrow\infty} \frac{\|q(t)\|_D}{t}  \geq \frac{C}{1+\varepsilon_0^{\zeta}}, \quad
 \limsup_{t\rightarrow\infty} \frac{\|q(t)\|_D}{t}  \leq \frac{C}{1-\varepsilon_0^{\zeta}}
$$
for some numerical constant $0<\zeta<1$.
\end{Theorem}

There are also some recent works on the ballistic motion for the discrete Schr\"{o}dinger operator on $\ell^2(\Z)$ with purely absolutely continuous spectrum.
This is described by the linear growth of the diffusion norm $\left(\sum_{n\in\Z}n^2| q_n(t)|^2\right)^{\frac12}$. See \cite{DamanikLY, Fillman, Kachkovskiy, ZhangZhao, Zhao} for details.

\

\noindent {\bf Idea of proof.} For the case where there is no potential, we can see the linear growth of diffusion norm simply by the Fourier transform(see Appendix \ref{fourier}).
In the presence of potential, our principal strategy is to relate the growth to the spectral transformation, instead of the Fourier transform.
Roughly speaking, for $g(E,t)=\int_{\R}q(x,t)\psi(E,x) dx$, with $\psi(E,x)$, $E\in\sigma(H)$, a generalized eigenvector of $H$, we have ${\rm i}\partial_t g(E,t) = E g(E,t)$, then
$$\int_{\R}q(x,t)\psi(E,x) dx=g(E,t) = e^{-{\rm i}E t}g(E,0).$$
So if $\psi(E,x)$ has nice differentiability and the derivative is well estimated, we can get
$$\int_{\R}q(x,t)\partial_E\psi(E,x) dx=\partial_E g(E,t)\sim t.$$
If, with some suitable measure $d\varphi$ supported on $\sigma(H)$, we have
\begin{equation}\label{intro_d_norm}
\left\|\int_{\R}q(x,t)\partial_E\psi(E,x) dx\right\|_{L^2(d\varphi)}\sim\left( \int_{\R}x^2|q(x,t)|^2 dx\right)^{\frac12},
\end{equation}
the linear growth of $\|q(t)\|_D$ is shown.

The above process can be realized by the ``modified spectral transformation", by generalizing the method in \cite{Zhao}. The generalized eigenvectors, with the Bloch-wave structures, are constructed by the previous works of Eliasson\cite{E92} for the reducibility of Schr\"odinger cocycle.
By adding some smoothing factors to the generalized eigenvectors(in a small part of the spectrum), the differentiability is improved.

The main difference from \cite{Zhao} is that the spectrum of linear operator is unbounded, so is the rotation number of the corresponding quasi-periodic cocycle.
To overcome this disadvantage, we construct a measure supported on the spectrum, with a factor decaying with respect to the spectral parameter $E$.

%This procedure makes the integration by parts over $[\inf\sigma(H),\infty)$ more complicated and more calculations need to be done.
%We will list several integrations in Appendix \ref{appendix_integral}, which are necessary in the proof.

In this way, the $L^2-$norm of the derivative(w.r.t. $E$) of the modified spectral transformation is close to the diffusion norm and (\ref{intro_d_norm}) is obtained.

\

The remaining part of paper is organized as follows.
Based on some basic notions on the Schr\"{o}dinger operator and Schr\"{o}dinger cocycle in Section \ref{sec_pre},
we introduce the previous works on the reducibility of Schr\"{o}dinger cocycle and some further properties in Section \ref{spectral_properties}.
This is devoted to compute some special integrals on an unbounded interval in Section \ref{sec_integration}.
Then we shall define the modified spectral transformation in Section \ref{section_spec_trans}, and finish the proof of Theorem \ref{thm-qp} in Section \ref{proof_theorem}.

\section{Preliminaries}\label{sec_pre}

\subsection{Schr\"{o}dinger operator and Schr\"{o}dinger cocycle}\label{pre_op_cocycle}
\noindent

In this subsection, we recall some basic notions and results for the quasi-periodic Schr\"{o}dinger operator
$H:L^2(\R) \rightarrow L^2(\R)$,
$$
(Hq)(x)=-q''(x)+ V(\omega x) q(x), \quad x\in\R,
$$
with $V$ and $\omega$ given as in (\ref{QP_Sch_EQ}).
Since here the potential is bounded, $H$ is essentially self-adjoint on ${\cal C}_c^{\infty}(\R)$, i.e., the subspace of smooth functions with compact support, and we can interpret $H$ as the unique extension to $L^2(\R)$.
It is known that its spectrum $\sigma(H)\subset [0,\infty)$ is unbounded.

We also introduce the corresponding Schr\"odinger cocycle $(\omega, A_0+F_0)$:
$$
\left(\begin{array}{c}
q \\
q'
\end{array}
\right)'=\left(A_0(E)+F_0(\omega x)\right)\left(\begin{array}{c}
q \\
q'
\end{array}
\right),
$$
with $A_0(E):=\left(\begin{array}{cc}
            0 & 1 \\
            -E & 0
          \end{array}
\right)$ and $F_0(\omega x):=\left(\begin{array}{cc}
            0& 0 \\
            V(\omega x) & 0
          \end{array}
\right)$.
Note that $(\omega, A_0+F_0)$ is equivalent to the eigenvalue problem $Hq=Eq$.

%\subsubsection{Spectral measure and integrated density of states}
%\noindent
%
%For any $\psi\in L^2(\R)$, let $\mu=\mu_{\psi}$ be
%the spectral measure of $H$ corresponding to $\psi$, which is defined so that
%$$\la(H-E)^{-1}\psi, \psi \ra = \int_\R\frac{1}{E-E'}d\mu_{\psi}(E'), \quad \forall \, E\in \C\setminus\sigma(H).$$

\subsubsection{Rotation number of Schr\"{o}dinger cocycle}
\noindent

We follow the same presentation as in \cite{HA2012} to give the precise definition of the rotation number for the time-continuous Schr\"{o}dinger cocycle $(\omega, A_0+F_0)$.

Let $\phi_E:\R\rightarrow SL(2,\R)$ be continuous with $\phi_E(0)=Id.$.
Given any $X\in \R^2\setminus\{0\}$, let
$\varphi_E(t,X)=\arg(\phi_E(t)X)$,
i.e.,
$$\varphi_E(t_2,X)-\varphi_E(t_1,X)=\Im\left(\int_\Gamma\frac{1}{z}dz\right),$$
where $\Gamma$ is the curve $\phi_E(t)X$, $t_1 \leq t \leq t_2$, in the complex plane.
We fix $0\leq \arg(X)<2\pi$ so that $\varphi_E$ is a continuous single-valued function in $\R\times(\R^2\setminus([0,\infty)\times\{0\}))$.

Suppose that $\varphi_E(t,X)-\varphi_E(t,Y)=2n\pi$.
This implies that $\phi_E(t)X=a \phi_E(t)Y$ for some $a>0$, i.e., $\phi_E(t)(X-aY)=0$. Hence $X=aY$. It follows that
\begin{equation}\label{Eq1.4_sana}
|\varphi_E(t,X)-\varphi_E(t,Y)|<\pi, \quad \forall \,  X,Y\in\R^2\setminus\{0\}.
\end{equation}
We say that the cocycle $(\omega, \phi_E)$ has the rotation number $\rho(E)$ if
$$\lim_{t\rightarrow\infty}\frac{\varphi_E(t,X)}{t}=\rho(E)$$
for some, and hence for all, $X \in \R^2\setminus\{0\}$. (\ref{Eq1.4_sana}) implies that the convergence is uniform with respect to $X$. If $(\omega, \phi)$ has the rotation number $\rho$, then
$$\lim_{t\rightarrow\infty}\frac{\varphi_E(t, X) - \varphi_E(t_0, X)}{t}= \rho(E),$$
which shows that the rotation number does not depend on the values of $\phi_E(t)$ over a finite interval.
These limits exist and define a continuous function of $E\in\R$.
For more details, we can refer to \cite{E92, JM}.

For the Schr\"{o}dinger cocycle $(\omega, A_0+F_0)$, the gap-labelling theorem relates the rotation number to the spectral property of the corresponding Schr\"{o}dinger operator:
\begin{Theorem}\cite{JM}\label{gap-labelling}
 Given any open interval $\cal I$ contained in $[\inf\sigma(H),\infty)\setminus\sigma(H)$, there is a unique $l\in\Z^d\setminus\{0\}$ such that $\rho=\frac{\la l, \,  \omega\ra}{2}$ on $\cal I$.
\end{Theorem}

\subsubsection{Classical spectral transformation}
\noindent

%Let $H:\ell^2(\Z)\rightarrow \ell^2(\Z)$ be the one-dimensional QP Schr\"odinger operator defined as in (\ref{qpSch_op}).
%Given $z\in\C$ with $\Im z > 0$ and fixed $\theta\in\T^d$, it is shown that the difference equation
%\begin{equation}\label{evproblem}
%-(q_{n+1}+q_{n-1})+V(\theta+n\omega)q_n=zq_n
%\end{equation}
%has unique solutions $q^{\pm}(z)$ which are $\ell^2$ at $\pm\infty$ respectively.
%Then the $m-$functions are defined as
%$$m^+(z)=-\frac{q_1^+(z)}{q_0^+(z)}, \quad m^-(z)=\frac{q_1^-(z)}{q_0^-(z)},$$
%which introduce the spectral measures $\mu^+$ and $\mu^-$ for the half-line Schr\"odinger operators on $\ell^2(\Z_+)$ and $\ell^2(\Z_-\cup\{0\})$ respectively by
%$$m^+(z)=\int_{\R}\frac{d\mu^+(E)}{E-z},\quad m^-(z)=\int_{\R}\frac{d\mu^-(E)}{E-z}.$$
%

Let $u(E,x)$ and $v(E,x)$ be solutions of the eigenvalue problem $Hq=Eq$ such that
\begin{equation}\label{initial_uv}
\left(
\begin{array}{cc}
u(E,0) & v(E,0)  \\
\partial_xu(E,0)  & \partial_x v(E,0)  \\
\end{array}
\right)=\left(
\begin{array}{cc}
1 & 0 \\
0 & 1 \\
\end{array}
\right).
\end{equation}

\begin{Theorem}[Chapter 9 of \cite{CL}]\label{spectral_measure_matrix}
There exists a non-decreasing Hermitian matrix $\mu=(\mu_{jk})_{j,k=1,2}$ whose elements are of bounded variation on every bounded interval on $\R$, %satisfying
%$$\mu_{jk}(E_2)-\mu_{jk}(E_1)=\lim_{\epsilon\rightarrow 0_+}\frac{1}{\pi}\int_{E_1}^{E_2}\Im M_{jk}(\nu+i\epsilon)d\nu,$$
%at points of continuity $E_1$, $E_2$ of $\mu_{jk}$, where on $\H$,
%$$M=\left(\begin{array}{cc}
%            M_{11} & M_{12} \\
%            M_{21} & M_{22}
%          \end{array}
%\right):=-\frac{1}{m^{+}+m^{-}}\left(\begin{array}{cc}
%                                  1 & m^+ \\
%                                  -m^- & -m^+m^-
%                                \end{array}
%\right),$$
%$$M_{11}:=\frac{1}{m_{+}+m_{-}},\quad M_{12}=M_{21}:=\frac12\, \frac{m_{-}+m_{+}}{m_{-}-m_{+}},\quad M_{22}:=\frac{m_{-}\, m_{+}}{m_{-}-m_{+}},$$
such that for any $q\in L^2(\R)$, with
$(g_1(E),g_2(E)):=\left(\int_{\R} q(x) u(E,x) dx, \int_{\R} q(x) v(E,x) dx\right)$,
we have Parseval's equality
$$\int_{\R}|q(x)|^2 dx=\int_{\R}\sum_{j,k=1}^2 g_j(E)\bar g_k(E) d\mu_{jk}(E).$$
\end{Theorem}

Given any matrix of measures on $\R$ $d\varphi=\left(\begin{array}{cc}
                 d\varphi_{11} & d\varphi_{12} \\[1mm]
                 d\varphi_{21} & d\varphi_{22}
               \end{array}\right)$,
let ${\cal L}^2(d\varphi)$ be the space of vectors
$G=(g_j)_{j=1,2}$, with $g_j$ functions of $E\in\R$ satisfying
\begin{equation}\label{gen_L2}
\|G\|_{{\cal L}^2(d\varphi)}^2:=\sum_{j,k=1}^2 \int_\R g_j \, \bar g_k \,  d\varphi_{jk}<\infty.
\end{equation}
In view of Theorem \ref{spectral_measure_matrix}, the map
$q\mapsto \left(\begin{array}{c}
                                                         \int_{\R}q(x) u(E,x)dx \\[1mm]
                                                         \int_{\R}q(x) v(E,x)dx
                                                       \end{array}
     \right)$
defines a unitary transformation between $L^2(\R)$ and ${\cal L}^2(d\mu)$.
We call it as the {\bf classical spectral transformation}.

%By Chapter $\uppercase\expandafter{\romannumeral5}$ of \cite{PF}(Page 297), we know that
The matrix of measures $(d\mu_{jk})_{j,k=1,2}$ is constructed via $m-$functions. It is Hermitian-positive, and
therefore each $d\mu_{jk}$ is absolutely continuous with respect to the measure $d\mu_{11}+d\mu_{22}$.
This measure determines the spectral type of the operator.
In particular, if the spectrum of $H$ is purely absolutely continuous, we have, for any $q\in L^2(\R)\setminus\{0\}$,
the classical spectral transformation is supported on a subset of $\sigma(H)$ with positive Lebesgue measure.

\

As mentioned in the previous section, to see the growth of the diffusion norm, we need the differentiability with respect to $E$.
But for the classical spectral transformation, there are some singularities with respect to $E$. More precisely, $u(x,E)$ and $v(x,E)$ are not well differentiated with respect to $E$ somewhere in the spectrum $\sigma(H)$.

For example, if we consider the Laplacian
$(H q)(x)=-q''(x)$, which may be the simplest case,
the linear growth of the diffusion norm of the solution to the equation ${\rm i} \partial_tq(x,t)=-\partial_x^2q(x,t)$ can be verified by the Fourier transform(see Appendix \ref{fourier}).
On the other hand, we have $\sigma(H)=[0,\infty)$ and for $E\in\sigma(H)$ the rotation number is $$\xi_0(E)=\rho_{(\omega, A_0)}(E)=\sqrt{E}\in[0,\infty).$$
It is easy to verify that the two generalized eigenvectors
\begin{equation}\label{classical_uv}
u(E,x)=\cos(\sqrt{E}x), \quad v(E,x)=\frac{\sin(\sqrt{E}x)}{\sqrt{E}}
\end{equation}
satisfy (\ref{initial_uv}). We can see that
the singularity comes when $\xi_0$ approaches $0$.

The matrix of spectral measures is $d\varphi=\frac{1}{2\pi}\left(
\begin{array}{cc}
\frac{1}{\sqrt{E}}dE & 0 \\
0 & \sqrt{E} dE \\
\end{array}
\right)$ on $(0,\infty)$(see Example 1, Page 252 of \cite{CL}).
Then for any $q\in L^2(\R)$, we have
\begin{eqnarray}
\left\|\left(\begin{array}{c}
\int_{\R}q(x) u(E,x)dx \\[1mm]
\int_{\R}q(x) v(E,x)dx
\end{array}
\right)\right\|^2_{{\cal L}^2(d\varphi)}
&=&\frac{1}{2\pi}\left[\int_{\R^2} q(x)\bar q(y) \int_0^\infty \cos(\sqrt{E}x) \cos(\sqrt{E}y) dx\, dy  \frac{dE}{\sqrt{E}}\right. \nonumber\\
& & + \, \left. \int_{\R^2} q(x)\bar q(y) \int_0^\infty \frac{\sin(\sqrt{E}x)}{\sqrt{E}} \frac{\sin(\sqrt{E}y)}{\sqrt{E}} dx\, dy  \sqrt{E} \, dE\right]\nonumber\\
&=& \frac{1}{2\pi} \int_{\R^2} q(x)\bar q(y) \int_0^\infty \cos(\sqrt{E}(x-y))  dx\, dy  \frac{dE}{\sqrt{E}}\nonumber\\
&=&\frac{1}{\pi} \int_{\R^2} q(x)\bar q(y) \int_0^\infty \cos(x-y)\rho \,  dx\, dy \,  \rho' dE \nonumber\\
&=&\int_{\R} |q(x)|^2 dx,\label{parseval_classic}
\end{eqnarray}
since it is well known that
$$\lim_{M\rightarrow\infty} \frac{1}{\pi} \int_{0}^M d\rho  \int_\R f(y) \cos(x-y)\rho\, dy =\frac12(f(x-0)+f(x+0))$$
for any a function of bounded variation in a neighbourhood of $x$.

\subsection{Regularity in the sense of Whitney}
\noindent

Given a closed subset $S$ of $\R$. We give a precise definition of ${\cal C}^1$ in the sense of Whitney, corresponding to a more general definition in \cite{Poschel}.

\begin{Definition}\label{definition_whitney}
Given two functions $F_0$, $F_1:S\rightarrow \C$(or $\M(2,\C)$) with some $0<M<\infty$, such that
$$
|F_0(x)|,\, |F_1(x)|\leq M, \;\  |F_0(x)-F_0(y)-F_1(y)(x-y)|\leq M|x-y|, \quad  \forall \, x,y\in S.
$$
We say that $F_0$ is {\bf ${\cal C}^1$ in the sense of Whitney} on $S$, denoted by $F_0\in{\cal C}_W^1(S)$, with the first order derivative $F_1$.
The ${\cal C}_W^1(S)-$norm of $F_0$ is defined as
$$|F_0|_{{\cal C}_W^1(S)}:=\inf M.$$
\end{Definition}

\begin{Remark}
By Whitney's extension theorem\cite{Whitney}, we can find an extension $\tilde F:\R\rightarrow \C$(or $\M(2,\C)$), which is ${\cal C}^1$ on $\R$ in the natural sense, such that
$\tilde F|_S=F_0$ and $\tilde F'|_S=F_1$.
\end{Remark}

\subsection{Notations}

\noindent {\rm 1)} With $\omega$ the Diophantine vector as above, we denote $\la k\ra:=\frac{\la k,\,  \omega\ra}{2}$ for any $k\in\Z^d$.

\smallskip

\noindent {\rm 2)} For any subset $S\subset\R$, let $\sharp(S)$ denote its cardinality of set, $\partial S$ be the set of its endpoints, $|S|$ be its Lebesgue measure, and $\rho(S)$ be its image by $\rho=\rho_{(\omega,\, A_0+F_0)}$.
%$S_1\triangle S_2:=(S_1\cup S_2)\setminus(S_1\cap S_2)$ is the difference of any two subsets $S_1, S_2\in\R$.
 \begin{itemize}
   \item Given any function $F$ on $S \times(2\T)^d$, possibly matrix-valued, let $$|F|_{S,\, (2\T)^d}:=\sup_{E\in S}\sup_{\theta\in(2\T)^d}|F(E,\theta)|.$$
If $F$ is ${\cal C}_W^1$ on $S$, then we define
    $|F|_{{\cal C}_W^1(S),\, (2\T)^d}:=\sup_{\theta\in(2\T)^d}|F(\cdot,\theta)|_{{\cal C}_W^1(S)}$.
   \item If $F$ is left and right continuous on $E$, then
   $F(E\pm):=\lim_{\epsilon\rightarrow0+}F(E\pm\epsilon)$.
   On the interval $(E_1, E_2)\subset\R$, if $F$ is left and right continuous on $E_1$ and $E_2$, then
   $$\left. F\right|_{[E_1, E_2]}=\left. F\right|_{(E_1, E_2)}:=F(E_{2}-)-F(E_1+),\; \left. F\right|^{E^+_{2}}_{E^-_1}:=F(E_{2}+)-F(E_1-).$$
  % In particular, if $F$ is continuous on $E_1$ and $E_2$, then $\left. F\right|^{E_{2}}_{E_1}:= F(E_{2})- F(E_1)$.
 \end{itemize}

\smallskip

\noindent {\rm 3)} For the quantities depending on the spectral parameter $E\in\R$, we do not always present this dependence explicitly and we simplify the notation ``$\partial_E$" into $\partial$, which denotes the derivative in the sense of Whitney on a certain subset of $\R$.

\section{Reducibility of quasi-periodic Schr\"odinger cocycle}\label{spectral_properties}
\noindent

Based on the general notions of Schr\"odinger operator and Schr\"odinger cocycle given in the previous section,
we present some further spectral properties, under the assumption that the potential function $V$ is sufficiently small.

\subsection{KAM scheme for the reducibility}
\noindent

In this subsection, we review the KAM theory of Eliasson\cite{E92} for the reducibility of the cocycle $(\omega, \, A_0+F_0)$,
which improves the previous results of Dinaburg-Sinai\cite{DiSi} and Moser-P\"oschel\cite{MP}.

With $\varepsilon_0=|V|_r$, $\sigma=\frac{1}{50}$, we define the sequences as in \cite{E92}:
$$\varepsilon_{j+1}=\varepsilon_{j}^{1+\sigma}, \quad N_j=4^{j+1}\sigma|\ln\varepsilon_j|, \quad j\geq 0.$$
Corresponding to the cocycle $(\omega, A_0+F_0)$, the result is formulated by the matrix equation
\begin{equation}\label{qpcocycle}
X'(x)=(A_0(E)+F_0(\omega x))X(x),
\end{equation}
recalling that $A_0(E):=\left(\begin{array}{cc}
            0 & 1 \\
            -E & 0
          \end{array}
\right)$ and $F_0(\omega x):=\left(\begin{array}{cc}
            0& 0 \\
            V(\omega x) & 0
          \end{array}
\right)$.

\begin{Proposition} \label{propHakan} There exists $\varepsilon_*=\varepsilon_*(\gamma,\tau, r)$ such that if $|V|_r=\varepsilon_0\leq\varepsilon_*$,
 then there is a full-measure subset $\Sigma=\cup_{j\geq0}\Sigma_j$ of $\sigma(H)$ with $\{\Sigma_{j}\}_j$ mutually disjoint Borel sets, satisfying
 \begin{equation}\label{measure_sigma_j}
 |\rho\left(\Sigma_{j+1}\right)|\leq |\ln\varepsilon_0|^{(j+1)^3 d} \varepsilon_{j}^{\sigma}, \quad  j\geq 0,
 \end{equation}
such that the following statements hold.
 \begin{itemize}
\item [(1)] The Schr\"odinger cocycle $(\omega, A_0+F_0)$ is {\bf reducible} on $\Sigma$. More precisely,
there exist
$\left\{ \begin{array}{l}
            B:\Sigma\rightarrow sl(2,\R) \;\  with \;\ eigenvalues \;\  \pm{\rm i}\rho\\[1mm]
            Y:\Sigma\times(2\T)^d\rightarrow GL(2,\R) \;\ analytic \;\  on \;\ (2\T)^d
          \end{array}
 \right.$
such that $\displaystyle X(x)=Y(\omega x) e^{Bx}$ is a solution of (\ref{qpcocycle}).
\item [(2)] For every $j\geq 0$, there is $k_j:\Sigma\rightarrow\Z^d$, such that
\begin{itemize}
  \item $|k_l|_{\Sigma_j}=0$ if $l\geq j$,
  \item $0<|k_j|\leq N_{j}$ on $\Sigma_{j+1}$ and $0<|\rho-\sum_{l\geq0}\la k_l\ra|_{\Sigma_{j+1}}< 2 \varepsilon_{j}^{\sigma}$.
\end{itemize}
\item [(3)] $B$ and $Y$ are ${\cal C}^1_W$ on ${\Sigma}_0$, and, with $\xi:=\rho-\sum_{j\geq0}\la k_j\ra$, $s\geq 2$, $\xi^{s+2} B$ and $\xi^{s+2} Y$ are ${\cal C}^1_W$
 on each ${\Sigma}_{j+1}$, $j\geq 0$. Moreover,
\begin{equation}\label{limit_state_whitney}
\left\{\begin{array}{lll}
|Y-Id.|_{{\cal C}^1_W({\Sigma}_0),\, (2\T)^d}, & |B-A_0|_{{\cal C}^1_W({\Sigma}_0)}\leq \varepsilon_0^{\frac13} & \\[1mm]
|\xi^{s+2\nu}Y|_{{\cal C}^\nu_W({\Sigma}_{j+1}),\, (2\T)^d}, & |\xi^{s+2\nu} B|_{{\cal C}^\nu_W({\Sigma}_{j+1})}\leq \varepsilon_j^{\frac{2\sigma}{3}}, & \nu=0,1
\end{array}\right. .
\end{equation}
%\item [(4)] For every $k\in\Z^d\setminus\{0\}$, there are two positive numerical constants $c'$, $\iota$, and $\beta'=\beta'(\gamma,\tau,r)$ such that
%\begin{equation}\label{gap-estime}
%|\rho^{-1}(\la k \ra)|\leq c' e^{-\beta'|k|^\iota}.
%\end{equation}
 \end{itemize}
\end{Proposition}

%We can refer to Remark 3.1--3.3 in \cite{Zhao} for more descriptions about this KAM regime and the (almost) reducibility for the Schr\"odinger cocycle, even though they are written for the $SL(2,\R)$ case.

\begin{Remark}
The original form Eliasson's theorem is: the cocycle $(\omega,\, A_0+F_0)$ is reducible if the rotation number $\rho$ is Diophantine
or rational with respect to $\frac{\omega}{2}$.
\begin{itemize}
  \item ``Rational w.r.t. $\frac{\omega}{2}$" means $\rho=\la k \ra$ for some $k\in\Z^d$. By the gap-labelling theorem, this case corresponds to the energies in $\R\setminus\sigma(H)$, where we have the reducibility to a hyperbolic matrix.
Hence, the corresponding eigenvectors constructed via the cocycle are exponentially growing/decaying vectors.
  \item In contrast, ``Diophantine w.r.t. $\frac{\omega}{2}$" means there exist $\gamma$, $\tau>0$ such that $|\rho-\la l\ra|>\frac{\gamma}{|l|^\tau}$ for any $l\in\Z^d\setminus\{0\}$. This corresponds to the energies in a full-measure subset of $\sigma(H)$.
The corresponding eigenvectors are the ``Bloch-waves"(see Subsectin \ref{subsec_bloch}).
\end{itemize}
\end{Remark}

\begin{Remark}\label{rmq_Prop1}
Associated with the above Diophantine condition, if, in $\sigma(H)$, $\rho$ is well separated from $\{\la l\ra\}_{l\in\Z^d\setminus\{0\}}$, it is the idealest case for applying the KAM scheme.
\begin{itemize}
  \item On $\Sigma_0$, there is no resonance for the rotation number $\rho$, so the standard KAM iteration is always applicable.
 $\Sigma_0$ is exactly the positive-measure subset of parameters for reducibility in the result of Dinaburg-Sinai\cite{DiSi}.
  \item On $\Sigma_{j+1}$, $j\geq 0$, there is some $k\in\Z^d$ with $0<|k|\leq N_{j+1}$, which appears as $k=\sum_{l=0}^j k_l$, such that
$0<|\rho-\la k\ra|_{\Sigma_{j+1}}< 2 \varepsilon_{j}^{\sigma}$.
But the resonance stops exactly at the $j^{\rm th}-$KAM step.
We could also apply the standard KAM on these subsets from the $(j+1)^{\rm th}-$step, because we could renormalize $\rho$ into $\xi:=\rho-\la k\ra$(the renormalization is done by several steps), which is well separated from $\{\la l\ra\}_{l\in\Z^d\setminus\{0\}}$.
Note that the ``renormalized rotation number" $\xi$ is close to $0$ on $\Sigma_{j+1}$ and it vanishes on the spectral gap where $\rho=\la k\ra$. So it can serve as a ``smoothing factor" on $\Sigma_{j+1}$.
\end{itemize}
Because of the difference between the procedures on $\Sigma_0$ and $\Sigma_{j+1}$, the transformation $Z$ and the reduced matrix $B$ possess different properties. In particular, on $\Sigma_{j+1}$, there are singularities like $\sim\xi^{-1}$(and $\sim \xi^{-3}$ after the derivation) for $Z$ and $B$.
Then, by multiplying $\xi^{s}$, $s\geq4$ the regularity is well improved as in (\ref{limit_state_whitney}).
Indeed, to get the ${\cal C}^1_W$ regularity, $\xi^3$ is enough, and the $4^{\rm th}$ power makes the norms small.
For better regularity, higher power of $\xi$ is needed.
\end{Remark}

%\begin{Remark}
%It has been shown in \cite{E92} and \cite{HA} that, for any $E\in\sigma(H)$, the Schr\"odinger cocycle $(\omega, A_0+F_0)$ is {\bf almost reducible}, i.e.,
%we can transform it arbitrarily close to a constant cocycle by a sequence of conjugations, without verifying the convergence of this sequence. On $\Sigma_j$, $j\geq 0$, since the resonance stops at exactly the $j^{\rm th}-$step and afterwards the conjugations are all close to identity, the convergence of sequence of conjugations is shown. Hence, in particular, reducibility holds for a.e. $E\in\sigma(H)$.
%\end{Remark}

\

Given $M\in\R$ satisfying $|M|>1$, with $J=J(M):=\min\left\{j\in\N: \; |M|\leq\varepsilon_{j}^{-\sigma}\right\}$, an approximation for the reducibility of Schr\"odinger cocycle $(\omega, A_0+F_0)$ can be stated in the following way, which will be contributed to computing several integrals on $[\inf\sigma(H), \infty)$(see Section \ref{sec_integration}).
\begin{Proposition}\label{propHakan1}
Let $|V|_r=\varepsilon_0\leq \varepsilon_*$ be as in Proposition \ref{propHakan}.
There is $$\Gamma^{(M)}=\bigcup_{j=0}^{J+1}\Gamma^{(M)}_j\subset[\inf\sigma(H), \infty),$$
with $\{\Gamma^{(M)}_j\}_{j=0}^{J+1}$ mutually disjoint and $\Sigma_{j}\subset \Gamma^{(M)}_{j}$, satisfying
\begin{equation}\label{measure}
\sharp\left([\inf\sigma(H), \infty)\setminus \Gamma^{(M)}\right)\leq |\ln\varepsilon_0|^{(J+1)^3 d}, \;\ \left|\rho\left(\Gamma^{(M)}_{j+1}\right)\right|\leq |\ln\varepsilon_0|^{(j+1)^3 d} \varepsilon_{j}^{\sigma},
\end{equation}
and $\left\{ \begin{array}{l}
                              A^{(M)}:\Gamma^{(M)}\rightarrow sl(2,\R) \;\ with \;\ two \;\ eigenvalues \;\ \pm{\rm i}\alpha^{(M)} \\[1mm]
                              Y^{(M)}:\Gamma^{(M)}\times(2\T)^d\rightarrow GL(2,\R) \;\ analytic \;\  on \;\ (2\T)^d
                            \end{array}
 \right.$, such that the following statements hold.\\
\noindent {\bf (S1)} $|\Re\alpha^{(M)}-\rho|_{\Gamma^{(M)}}\leq \varepsilon_J^{\frac14}$ and for $0\leq j\leq J$, there is $k^{(M)}_j:\Gamma^{(M)}\rightarrow\Z^d$, constant on each connected component of $\Gamma^{(M)}$, such that
\begin{itemize}
\item [1.] $|k^{(M)}_{l}|_{\Gamma^{(M)}_{j}}=0$ if $l\geq j$,
\item [2.] $0<|k^{(M)}_j|\leq N_{j}$ on $\Gamma^{(M)}_{j+1}$ and $|\Re\alpha^{(M)} - \sum_{l=0}^{J}\la k^{(M)}_l\ra|_{\Gamma^{(M)}_{j+1}}\leq \frac32 \varepsilon_j^{\sigma}$, $0\leq j\leq J$.
\end{itemize}
{\bf (S2)} Let $\xi^{(M)}:=\Re\alpha^{(M)} - \sum_{l=0}^{J}\la k^{(M)}_l\ra$.
\begin{itemize}
\item On $\Gamma^{(M)}_{j+1}$, $0\leq j\leq J$, in each connected component,
      there is one and only one subinterval ${\cal I}$ such that
       $\xi^{(M)}=0$ on ${\cal I}$, and outside ${\cal I}$, $\xi^{(M)}\neq0$ with
    \begin{equation}\label{esti_plat}
    \frac13<\partial\xi^{(M)}\leq N_j^{4\tau} |\xi^{(M)}|^{-1}, \quad |\partial^2 \xi^{(M)}|\leq N_j^{8\tau}|\xi^{(M)}|^{-3}.
    \end{equation}
%    Moreover, $0\leq|{\cal I}|\leq \varepsilon_{j}^{2\sigma}$  \footnote{$|{\cal I}|=0$ means that this subinterval degenerates into one point.}.
\item On $\Gamma^{(M)}_{0}$, if $\xi^{(M)}\neq 0$, we have $\partial\xi^{(M)}= \frac{\partial \det A^{(M)}}{2\xi^{(M)}}> \frac1{3\xi^{(M)}}$. \footnote{Indeed, the only possibility that $\xi^{(M)}=0$ on $\Gamma_0$ is in the interval containing $\inf\sigma(H)$.}
\end{itemize}
{\bf (S3)} $|Y^{(M)}-Y|_{\Sigma_0,\, (2\T)^d}$, $|A^{(M)}- B|_{\Sigma_0} \leq \varepsilon_{J}^{\frac14}$, and for $0\leq j\leq J$,
\begin{equation}\label{error_J+1}
|\xi^{(M)}\, Y^{(M)}- \xi\, Y|_{\Sigma_{j+1},\, (2\T)^d},\quad |\xi^{(M)}\, A^{(M)}-\xi\, B|_{\Sigma_{j+1}}\leq \varepsilon_{J}^{\frac14}.
\end{equation}
For $\nu=0,1,2$,
\begin{equation}\label{sigma_m_0}
\left\{
\begin{array}{l}
\displaystyle |\partial^\nu (Y^{(M)}-Id.)|_{\Gamma^{(M)}_{0}, \, (2\T)^d},\; |\partial^\nu (A^{(M)}- A_0)|_{\Gamma^{(M)}_{0}}\leq \varepsilon_0^{\frac13}  \\[1mm]
\displaystyle |\partial^\nu Y^{(M)}|_{(2\T)^d}, \; |\partial^\nu A^{(M)}|\leq \frac{\varepsilon_j^{-\frac{\sigma}5}}{|\xi^{(M)}|^{1+2\nu}} \  on \ \Gamma^{(M)}_{j+1}  \  if  \  \xi^{(M)}\neq 0
\end{array}
\right. .
\end{equation}
{\bf (S4)} $\{E\in\partial\Gamma^{(M)}:M\rho(E)\notin\pi\Z\}\subset \partial\Gamma^{(M)}_{J+1}$.
     For any connected component $(E_*, E_{**})$ of $\Gamma^{(M)}_{J+1}$, we have
     \begin{equation}\label{connec_comp_J+1}
      \left|\left. \rho\right|_{(E_*, E_{**})}\right|\leq 2 \varepsilon_{J}^{\sigma(1+\frac\sigma2)}, \quad \varepsilon_J^{3\sigma(1+\sigma)}\leq E_{**}-E_{*}\leq \varepsilon_J^{\sigma(1+\frac{\sigma}3)}.
     \end{equation}
     Moreover, $k_j^{(M)}(E^-_*)=k_j^{(M)}(E^+_{**})$, $0\leq j\leq J$, and
     there is $0\leq j_* < J$ such that $E_*$, $E_{**}\in \partial\Gamma^{(M)}_{j_*}$, with
\begin{equation}\label{edge_point}
\left\{\begin{array}{llll}
\left|\left.(Y^{(M)}-Id.)\right|^{E^+_{**}}_{E^-_*}\right|_{(2\T)^d},   & \left|\left.( A^{(M)}-A_0)\right|^{E^+_{**}}_{E^-_*} \right|   &  \leq \displaystyle \frac{\varepsilon_0^{\frac13}}2(E_{**}-E_*), & j_*=0 \\[3mm]
\left|\left.(\xi^{(M)})^4\, Y^{(M)}\right|^{E^+_{**}}_{E^-_*} \right|_{(2\T)^d},   & \left|\left.(\xi^{(M)})^4\, A^{(M)}\right|^{E^+_{**}}_{E^-_*}\right|  &  \leq \displaystyle \frac{\varepsilon_{j_*-1}^{\frac{2\sigma}{3}}}2(E_{**}-E_*), & j_*\geq1
\end{array}\right.  .
\end{equation}
\end{Proposition}

%\begin{Remark}
%Since $\{k^{(M)}_j\}_{0\leq j\leq J}$ are piecewise constant, we have $\partial \rho^{(M)}=\partial\xi^{(M)}$,
%and equals to $-\frac{\partial{\rm tr}A^{(M)}}{2\sin\rho^{(M)}}$ if $\xi^{(M)}\neq 0$.
%\end{Remark}
%\begin{Proposition}
%Given any $m_1,\, m_2\in\Z$, with $l(m_1)\leq l(m_2)$.
%\begin{itemize}
%  \item [1)] For each $0\leq j\leq l(m_1)$, $\left|\Sigma_j^{(m_1)}\setminus\Sigma_j^{(m_2)}\right|\leq \varepsilon_{l(m_1)}^{\frac\sigma3}.$
%  \item [2)] One can find a subset ${\cal O}_{m_1,\, m_2}$ satisfying $\left|[\inf\sigma(H),\, \sup\sigma(H)]\setminus {\cal O}_{m_1,\, m_2}\right|\leq \varepsilon_{l(m_1)}^{\frac\sigma3}$ such that
%\begin{equation}\label{cauchyAZ}
%\left|\partial^\nu (\xi^{(m_1)}- \xi^{(m_2)})\right|,\;\ \left|\partial^\nu (\tilde A^{(m_1)}- \tilde A^{(m_2)})\right|,\;\ \left|\partial^\nu (\tilde Z^{(m_1)}(\theta)- \tilde Z^{(m_2)}(\theta))\right|_{\T^d}\leq \varepsilon_{l(m_1)}^{\frac\sigma3}.
%\end{equation}
%\end{itemize}
%\end{Proposition}

\begin{Remark}
The mutually disjoint subsets $\{\Gamma^{(M)}_j\}_{0\leq j\leq J+1}$ given in Proposition \ref{propHakan1} cover $[\inf\sigma(H), \infty)$ up to finite points(obviously, they are unions of finite intervals contained in $[\inf\sigma(H), \infty)$).
They divide the energies according to the extent of resonances.
As the iteration continues until the limit state, we can get the sequence of mutually disjointed subset $\{\Sigma_j\}_{j\geq 0}$ after excluding every gap in the spectrum.
\end{Remark}

\begin{Remark}
$Z^{(M)}$ and $A^{(M)}$ in Proposition \ref{propHakan1} are constructed by KAM iteration at the $J^{\rm th}-$step.
$A^{(M)}$ has two eigenvalues $\alpha^{(M)}$, with $\Re\alpha^{(M)}$ the approximation of the rotation number.
They are not uniquely determined, and in particular, as shown in {\bf (S4)}, for any given $M\in \R$ with $|M|>1$,
we can choose delicately the endpoints of the ``resonance intervals" at the initial several steps, such that $M\rho\in \pi\Z$ on these endpoints(since $\varepsilon_J^\sigma<\frac{1}{|M|}\leq\varepsilon^{\sigma}_{J-1}$, if $J>1$, the endpoints are adjustable within this range when $j< J$). We can refer to the proof of Proposition 1 in \cite{Zhao}.
\end{Remark}

\noindent
{\bf Modifications of proof for Proposition \ref{propHakan} and \ref{propHakan1}:}
The proofs are similar to that of Proposition 1 and 2 in \cite{Zhao},
since the original idea of the reducibility in \cite{HA} is similar to that of \cite{E92}.
There is a detailed proof for Proposition 1 and 2 of \cite{Zhao}(see Section 3.1 and Appendix A.1, A.2 of \cite{Zhao}), so we do not present it here, just give the modifications because of the difference between the continuous and the discrete Schr\"{o}dinger operators.
\begin{itemize}
  \item For the continuous Schr\"{o}dinger operator, the rotation number $\rho$ of the corresponding cocycle is unbounded, with the renormalized rotation number $\xi=\rho-\sum_{j\geq0}\la k_j\ra$ close to $0$ in a small part of the spectrum.
In this small subset, the singularity of the conjugation $Y$ and the constant matrix $B$ becomes ``$\sim\xi^{-1}$"(see (\ref{limit_state_whitney}) and (\ref{sigma_m_0})) instead of $\sim\sin^{-1}\xi$(see (3.1) and (3.24) in \cite{Zhao}).
  \item Because of the unboundedness of $\rho$, the transversality of the rotation number in the spectrum becomes weaker($\partial\rho>\frac1{2\rho}$ instead of $\partial\rho>\frac1{2\sin\rho}\geq \frac12$). It is similar for the approximated rotation number
on the large subset of $[\inf\sigma(H),\infty)$, i.e., $\Gamma_0^{(M)}$, as shown in {\bf (S2)} of Proposition \ref{propHakan1}.
Note that for the last connected component $(E_*, E_{**})\subset\Gamma^{(M)}_{0}$, we have $E_{**}=\infty$.

  \item In Proposition \ref{propHakan1}, we focus on the KAM iteration until the $(J+1)^{\rm th}-$step. So the connected components $(E_*, E_{**})\subset\Gamma_{J+1}^{(M)}$ are only related to the energies where the rotation number is close to $\la k\ra$, $|k|\leq N_J$. Hence $\partial\rho> \frac{c}{N_J}$ on $\sigma(H)\cap\Gamma_{J+1}^{(M)}$ with some constant $c>0$ depending on $\omega$. Combining with the $\frac12-$H\"older continuity of $\rho$ and the estimate of the spectral gaps(both of which will be given in the next subsection), we can still estimate the Lebesgue measure for each connected component of $\Gamma_{J+1}^{(M)}$ as in (\ref{connec_comp_J+1}).
\end{itemize}

%The proof of Proposition \ref{propHakan} and \ref{propHakan1} are based on the KAM mechanisms in \cite{E92, HA}.
%We leave the outline of their proofs in Appendix \ref{proof}.

\subsection{Further properties}\label{prop_rot_num}
\noindent

As a result of the reducibility, we have
\begin{Theorem}\cite{E92}
With $|V|_r=\varepsilon_0\leq \varepsilon_*$ as in Proposition \ref{propHakan}, the
spectrum of $H$ is purely absolutely continuous.
\end{Theorem}

For the rotation number $\rho=\rho_{(\omega,\, A_0+F_0)}$,
we also have the following further results, which come with the analysis on the reducibility of Schr\"odinger cocycle.
\begin{Theorem}\cite{HA2012}\label{proposition_rotation_number}
With $|V|_r=\varepsilon_0\leq \varepsilon_*$ as in Proposition \ref{propHakan}, we have
\begin{enumerate}
  \item $\rho=\rho_{(\omega,\, A_0+F_0)}$ is $\frac12-$H\"{o}lder continuous, i.e., there is a numerical constant $c>0$, such that for any given $E_1$, $E_2\in\R$,
  $$|\rho(E_1)-\rho(E_2)|<c |E_1-E_2|^{\frac12}.$$
  \item $\rho=\rho_{(\omega,\, A_0+F_0)}$ is absolutely continuous on $\R$, i.e., given finite intervals $\{{\cal I}_j\}_j$ on $\R$, for any $\eta>0$, there exists $\delta=\delta(\eta)>0$, such that if $\sum_j|{\cal I}_j| < \delta$ then $\sum_j \left|\left.\rho\right|_{{\cal I}_j}\right| <\eta$.
\end{enumerate}
\end{Theorem}
%\proof Recalling that $A_0(E)=\left(\begin{array}{cc}
%            -E & -1 \\
%            1 & 0
%          \end{array}
%\right)$ and $F_0(\theta)=\left(\begin{array}{cc}
%            V(\theta) & 0 \\
%            0 & 0
%          \end{array}
%\right)$,
%the H\"{o}lder continuity and absolute continuity are obtained as direct corollaries of Theorem 2 in \cite{HA} and Theorem 1 in \cite{HA2014} respectively. \qed

\begin{Theorem}\cite{DS, E92}\label{proposition_rotation_number}
With $|V|_r=\varepsilon_0\leq \varepsilon_*$ as in Proposition \ref{propHakan}, we have
\begin{equation}\label{transversality_rho}
(2\rho)^{-1}<\partial\rho<\infty \ for \ a.e. \ E\in\sigma(H).
\end{equation}
\end{Theorem}
%\proof According to Proposition \ref{spec_ac} in the previous subsection, if $|V|_r\leq \varepsilon_*$, then the spectrum of $H_\theta$ is purely absolutely continuous for any $\theta\in\T^d$.
%As the well-known result of Kotani theory, $L(E)=0$ for a.e. $E\in\sigma(H)$.
%In view of Theorem 1.4 of \cite{DS}, we get the conclusion.\qed

\begin{Theorem}\cite{HA2012}\label{gap-estimate}
With $|V|_r=\varepsilon_0\leq \varepsilon_*$ as in Proposition \ref{propHakan},
there are two positive numerical constants $c'$, $\iota$, and $\beta'=\beta'(\gamma,\tau,r)$ such that
\begin{equation}\label{gap-estime}
|\rho^{-1}(\la k \ra)|\leq c' e^{-\beta'|k|^\iota},\quad \forall \, k\in\Z^d.
\end{equation}
\end{Theorem}

In view of (\ref{measure}), we can estimate the Lebesgue measure of the subset $\Gamma_{j+1}^{(M)}$, $j\geq 0$, as $\left|\Gamma_{j+1}^{(M)}\right|< \varepsilon_0^{\frac\sigma2}$, by applying Theorem \ref{proposition_rotation_number} and \ref{gap-estimate}.

\

From now on, we always assume that $|V|_r=\varepsilon_0 <\varepsilon_*$, and $\varepsilon_0$ is small enough such that it is compatible with every simple calculation in this paper.
Moreover, we also assume that (\ref{transversality_rho}) is satisfied on the full-measure subset $\Sigma$ of $\sigma(H)$ given in Proposition \ref{propHakan}.

\subsection{Construction of Bloch-waves}\label{subsec_bloch}
\noindent

In general, the {\bf Bloch-wave} of a self-adjoint operator on $L^2(\R)$ means the generalized eigenvector $\psi$, of the form
$\psi(x)=e^{{\rm i}\varrho x} h(\tilde\alpha x)$, with $\varrho$, $\tilde\alpha$ some real numbers, and $h$ a periodic function of $x\in\R$.
Here $\varrho$ is usually called the {\bf Floquet exponent}.

Back to Proposition \ref{propHakan}, we can construct Bloch-waves of Schr\"odinger operator $H$ on $\Sigma$.
More precisely, for the Schr\"odinger operator $H$, by the matrices
$Y=\left(\begin{array}{cc}
     Y_{11} & Y_{12} \\
     Y_{21} & Y_{22}
   \end{array}\right)$ and
$B=\left(\begin{array}{cc}
B_{11} & B_{12} \\
B_{21} & B_{22}
\end{array}\right)$ given in Proposition \ref{propHakan},
we can see $\tilde\psi(E,x)=e^{{\rm i}x\rho(E)}\tilde f(E,\omega x)$ satisfies $H\psi=E\psi$, with $\tilde f:\Sigma\times (2\T)^d\rightarrow \C$ given by
$$\tilde f(E,\cdot):=Y_{11}\left(E,\cdot\right)B_{12}(E)+Y_{12}\left(E,\cdot\right)({\rm i}\rho(E)-B_{11}(E)).$$
Indeed, by noting that $\left(\begin{array}{c}
B_{12} \\[1mm]
{\rm i}\rho-B_{11}
\end{array}
\right)$ is an eigenvector of $B$ corresponding to the eigenvalue ${\rm i}\rho$,
we have
$$\left(\begin{array}{c}
\tilde\psi  \\[1mm]
\tilde\psi'
\end{array}
\right)=e^{{\rm i}x\rho}\, Y(\omega x)\left(\begin{array}{c}
B_{12} \\[1mm]
{\rm i}\rho-B_{11}
\end{array}
\right)=Y(\omega x)e^{Bx}\left(\begin{array}{c}
B_{12} \\[1mm]
{\rm i}\rho-B_{11}
\end{array}
\right).$$
Hence, we can also get the Bloch-wave
$$
\psi(x)=e^{{\rm i}x\rho}f(\omega x) \;\ {\rm with} \;\ f=\left\{\begin{array}{ll}
 \tilde f, & E\in\Sigma_0 \\[1mm]
\xi^8 \tilde f, & E\in\Sigma_{j+1} , \; j\geq0
\end{array} \right. .
$$

\begin{Remark} The Bloch-wave depends on the energy $E\in\Sigma$.
Recalling (\ref{limit_state_whitney}), we know that the Bloch-wave $\tilde\psi$ has nice estimates on $\Sigma_0$,
In contrast, it has some singularities ``$\sim\xi^{-1}$" on $\Sigma_{j+1}$, $j\geq 0$, whose union forms a small part of the spectrum.
So we add a smoothing factor $\xi^8$, just on this small part to cover the singularities and get better estimates.
\end{Remark}

Based on the Bloch-wave $\psi$, we can introduce the ingredients of the modified spectral transformation for the operator $H$(see Section \ref{section_spec_trans}).
Let $${\cal K}(x):={\Im}(\psi(x)), \;\  {\cal J}(x) :={\Re}(\psi(x)) \quad  {\rm for}  \ E\in\Sigma,$$
and ${\cal K}(x)|_{\R\setminus \Sigma}={\cal J}(x)|_{\R\setminus \Sigma}=0$. By a direct calculation, we find
$$
\psi(x)=e^{{\rm i}x\rho}f(\omega x) = \beta_0(x) e^{{\rm i}x\rho} +   {\rm i} \beta_1(x) \rho e^{{\rm i}x\rho},
$$
where, for $l=0, 1$, $\beta_{l}:\Sigma\times\R\rightarrow \R$, real-analytic on $\R$ and ${\cal C}^1_W$ on each $\Sigma_j$, $j\geq0$, is given by
$\beta_{l}=\left\{\begin{array}{ll}
\tilde \beta_{l}, & E\in\Sigma_0 \\[1mm]
     \xi^{8}\tilde \beta_{l},&  E\in\Sigma_{j+1},\;\ j\geq0
  \end{array}\right.$,
with
$$\tilde \beta_{0}(x) := Y_{11}(\omega x)B_{12}-Y_{12}(\omega x)B_{11}, \quad \tilde \beta_{1}(x) :=Y_{12}(\omega x).$$
So for $E\in\Sigma$, we have
$$
{\cal K}(x)= \beta_0(x)\sin (x\rho) +  \beta_1(x) \rho \cos (x\rho), \quad {\cal J}(x)= \beta_0(x) \cos (x\rho) -  \beta_1(x) \rho \sin (x\rho).
$$
According to (\ref{limit_state_whitney}) and the fact that $|\xi|\leq 2\varepsilon_j^{\sigma}$ on $\Sigma_{j+1}$, $j\geq 0$, it is obvious that
\begin{equation}\label{esti_abc_sigma_avant}
|\beta_{l}-\delta_{l,0}|_{{\cal C}^1_W(\Sigma_0),\, \R}\leq \varepsilon_0^{\frac14},\quad |\beta_{l}|_{{\cal C}^1_W(\Sigma_{j+1}),\, \R}\leq \varepsilon_j^{\sigma},\;\  j \geq 0.
%\left\{\begin{array}{ll}
%   |\beta_{n,n_\Delta}-\delta_{n,n_\Delta}|_{{\cal C}^1_W(\Sigma_0),\, (2\T)^d}\leq \varepsilon_0^{\frac14},  & E\in\Sigma_0 \\[1mm]
%   |\beta_{n,n_\Delta}|_{{\cal C}^1_W(\Sigma_{j+1}),\, (2\T)^d}\leq \varepsilon_j^{\sigma},  & E\in\Sigma_{j+1},\;\  j \geq 0
%\end{array} \right. .
\end{equation}
%Now, on $\Sigma$, we get two solutions of the eigenvalue problem $Hq=Eq$:
%\begin{equation}\label{u_v}
%u_n=\frac{{\Cal K}_n}{{\Cal K}_1}, \quad v_n=-\frac{{\cal J}_n}{{\Cal K}_1}, \quad n\in\Z,
%\end{equation}
%since they are both linear combinations of $\psi$ and $\bar \psi$.
%It can be verified that
%$$\left(\begin{array}{c}
%u_1 \\
%u_0
%\end{array}
%\right)=\left(\begin{array}{c}
%1 \\
%0
%\end{array}
%\right),\quad
%\left(\begin{array}{c}
%v_1 \\
%v_0
%\end{array}
%\right)=\left(\begin{array}{c}
%0 \\
%1
%\end{array}
%\right).
%$$

With $\Sigma_j$, $0\leq j \leq J+1$, $Y$, $B$, $\xi$ replaced by $\Gamma_j^{(M)}$, $Y^{(M)}$, $A^{(M)}$, $\xi^{(M)}$ given in Proposition \ref{propHakan1} respectively,
we get the piecewise ${\cal C}^2$ coefficients $\beta_{l}^{(M)}:\Gamma^{(M)}\times \R \rightarrow \R$, and
%
% get the approximations of ${\Cal K}_n$ and ${\cal J}_n$ as
%$${\Cal K}^{(M)}_n=\sum_{n_\Delta}\beta^{(M)}_{nn_\Delta} \sin n_\Delta\rho, \quad {\cal J}^{(M)}_n=\sum_{n_\Delta}\beta^{*(M)}_{nn_\Delta} \sin (n_\Delta-1)\rho,$$
%with the meaning of
\begin{Lemma}\label{coefficient_abc}
For $l=0, 1$,
$$\left\{ \begin{array}{ll}
            |\partial^\nu(\beta^{(M)}_{l} -\delta_{l, 0})|_{\Gamma^{(M)}_0,\, \R}\leq \varepsilon_{0}^{\frac14}, & \\[1mm]
            |\partial^\nu \beta_{l}^{(M)}|_{\R} \leq \varepsilon_j^{\frac\sigma5}|\xi^{(M)}|^{5-2\nu} \;\ on \;\ \Gamma^{(M)}_{j+1}, & 0\leq j\leq J
          \end{array}
 \right. ,\quad \nu=0,1,2,$$
% \begin{itemize}
%     \item On ${\cal I}^{(M)}\subset\Gamma^{(M)}_0$, $|\partial^\nu(\beta^{(M)}_{n n_\Delta} -\delta_{n_\Delta n})|_{\T^d}\leq \varepsilon_{0}^{\frac13}$.
%     \item On ${\cal I}^{(M)}\subset\Gamma^{(M)}_{j+1}$, $0\leq j\leq J$, $|\partial^\nu \beta_{n_\Delta}^{(M)}|_{\T^d} \leq \varepsilon_j^{\frac\sigma6}|\sin\xi^{(M)}|^{5-2\nu}$, $\nu=0,1,2$.
%   $$|\beta^{(M)}|\leq N_j^{3\tau}|\sin\xi^{(M)}|^6,\quad |\partial \beta^{(M)}|\leq 2\varepsilon_j^{-\frac\sigma5}|\sin\xi^{(M)}|^4,\quad |\partial^2 \beta^{(M)}|\leq 3\varepsilon_j^{-\frac\sigma2}|\sin\xi^{(M)}|^2.$$
%   \end{itemize}
and for each connected component $(E_*, E_{**})\subset \Gamma^{(M)}_{J+1}$,
 $ \left| \left.\beta^{(M)}_{l}\right|^{E^+_{**}}_{E^-_*} \right|_{\R} \leq \varepsilon_{J}^{\sigma(1+\frac\sigma4)}$.
\end{Lemma}

\proof We only prove for $\beta^{(M)}_{0}$, with that of $\beta^{(M)}_{1}$ similar.

On $\Gamma^{(M)}_0$, $\beta^{(M)}_{0}=\tilde \beta^{(M)}_{0}$ with
\begin{equation}\label{tilde_an_m}
\tilde \beta^{(M)}_{0}=Y^{(M)}_{11}\left(\omega x\right)A^{(M)}_{12}-Y^{(M)}_{12}\left(\omega x\right)A^{(M)}_{11}.
\end{equation}
Then, in view of (\ref{sigma_m_0}), $|\partial^\nu(\beta^{(M)}_{0} - 1)|_{\Gamma^{(M)}_0,\, \R}\leq \varepsilon_{0}^{\frac14}$ is evident.

On $\Gamma^{(M)}_{j+1}$, $0\leq j \leq J$, $\beta^{(M)}_{0}=(\xi^{(M)})^{8}\tilde \beta^{(M)}_{0}$.
In each connected component of $\Gamma^{(M)}_{j+1}$, as mentioned in {\bf (S2)}, there is a subinterval ${\cal I}$ where $\xi^{(M)}=0$.
So $\beta^{(M)}_{0}=0$ on ${\cal I}$.
Outside ${\cal I}$, $\xi^{(M)}\neq0$ then
by (\ref{esti_plat}), (\ref{sigma_m_0}) and (\ref{tilde_an_m}), we have, for $\nu=0,1,2$,
$$|\partial^{\nu}\tilde \beta^{(M)}_{0}|_{\R}\leq\varepsilon_j^{-\frac{2\sigma}{5}}|\xi^{(M)}|^{-(2+2\nu)},\quad
|\partial^{\nu} (\xi^{(M)})^{8}|\leq\varepsilon_j^{\frac{4}{5}\sigma}|\xi^{(M)}|^{7-2\nu}. $$
Hence, combining the estimates above, $|\partial^\nu \beta^{(M)}_{0}|_{\R} \leq \varepsilon_j^{\frac\sigma5}|\sin\xi^{(M)}|^{5-2\nu}$ on $\Gamma^{(M)}_{j+1}$.

\smallskip

For the connected component $(E_*, E_{**})\subset \Gamma^{(M)}_{J+1}$, according to {\bf (S3)}, there is $0\leq j_*\leq J$, such that $E^-_{*}, E^+_{**}\in \partial\Gamma^{(M)}_{j_*}$. By (\ref{sigma_m_0}) and (\ref{edge_point}),
\begin{itemize}
  \item if $j_*=0$,
$\left|\left.\beta^{(M)}_{0}\right|^{E^+_{**}}_{E^-_*}\right|_{\R}
=\left|\left.\tilde \beta^{(M)}_{0}\right|^{E^+_{**}}_{E^-_*} \right|_{\R}
\leq 10\, (E_{**}-E_*)\leq \varepsilon_J^{\sigma(1+\frac\sigma4)}$;
  \item if $j_* \geq 1$, then for $\beta^{(M)}_{0}=(\xi^{(M)})^{8}\tilde \beta^{(M)}_{0}$, $\left| \left.(\xi^{(M)})^{8}\tilde \beta^{(M)}_{0}\right|^{E^+_{**}}_{E^-_*} \right|_{\R}$ is bounded by terms like
$$4 \left| \left. (\xi^{(M)})^4\cdot Y^{(M)}\right|^{E^+_{**}}_{E^-_*} \right|_{\R}\cdot\left| (\xi^{(M)})^4\cdot A^{(M)} \right|_{\Gamma_{j_*}^{(M)},\, \R}\leq \varepsilon_J^{\sigma(1+\frac\sigma4)}.$$\qed
\end{itemize}

Moreover, by (\ref{error_J+1}), it is obvious that
\begin{equation}\label{error_beta}
|\beta_{l}-\beta^{(M)}_{l}|_{\Sigma_j,\, \R}\leq 10 \varepsilon_J^{\frac14}, \quad 0\leq j\leq J+1.
\end{equation}

%\begin{Lemma}\label{lemma_nn}
%For $l=0, 1$,
%$\displaystyle \left|\frac{8\varepsilon_0^{-\frac{7\sigma}8}}{\pi \csc(\frac{\pi}8)}\int_{\Sigma}\frac{\beta_l(x) \beta_l(y)}{1+\rho^8}\, \partial \rho\, dE- \delta_{l, 0}\right|\leq \varepsilon_0^{\sigma}$.
%\end{Lemma}
%\proof By (\ref{esti_abc_sigma_avant}), we can get for any $x, y \in\R$,
%$$\left|\beta_l(x) \beta_l(y)-\delta_{l, 0 }\right|_{\Sigma_0}\leq 3\varepsilon_0^{\frac14};\quad |\beta_l(x) \beta_l(y)|_{\Sigma_{j+1}}\leq \varepsilon_j^{2\sigma},\;\ j \geq 0.$$
%Then, recalling that $|\rho\left(\Sigma_{j+1}\right)|\leq |\ln\varepsilon_0|^{(j+1)^3 d} \, \varepsilon_{j}^{\sigma}$,
% we have $$\left|\int_{\Sigma_j} \frac{\beta_l(x) \beta_l(y)-\delta_{l, 0 }}{1+\rho^8} \, \partial \rho\,dE\right|\leq
%\left\{ \begin{array}{cc}
%          3\pi\varepsilon_0^{\frac14}, & j=0 \\[1mm]
%          \varepsilon_0^{\sigma}(1+\varepsilon_{j-1}^{2\sigma})\cdot|\ln\varepsilon_0|^{j^3 d} \, \varepsilon_{j-1}^{\sigma}, &  j\geq 1
%        \end{array}
%\right. .$$
%Therefore, by noting that $\displaystyle \int_{0}^\infty\frac{d\rho}{1+\rho^8}=\frac{\pi \csc\frac{\pi}{8}}{8}\varepsilon_0^{\frac{7\sigma}8}$, we have
%$$\left|\frac{8\varepsilon_0^{-\frac{7\sigma}8}}{\pi \csc(\frac{\pi}8)}\int_{\Sigma} \frac{\beta_l(x) \beta_l(y)}{1+\rho^8}\, \partial \rho\,dE-\delta_{l,0}\right|=\frac{8\varepsilon_0^{-\frac{7\sigma}8}}{\pi \csc(\frac{\pi}8)} \left| \int_{\Sigma} \frac{\beta_l(x) \beta_l(y)-\delta_{l, 0 }}{1+\rho^8} \, \partial \rho\,dE\right| \leq  \varepsilon_0^{\sigma}.$$
%\qed

\section{Some integrals on $[\inf\sigma(H),\infty)$}\label{sec_integration}
\noindent

In Proposition \ref{propHakan1}, we have divided the interval $[\inf\sigma(H), \infty)$ into $J(M)+2$ parts for some given $M\in\R$ with $|M|>1$, up to a subset of finite points.
With this division, we can estimate the following integrals, which will be applied in analyzing the modified spectral transformation in the next section.

For $l=0,1$, let $h_{l}=h_{l}(x,y,M):=\beta^{(M)}_l(x)\beta^{(M)}_l(y)$ for any $x,y\in\R$.
\begin{Proposition}\label{hmn_integral}
%Given $M\in\Z\setminus\{0\}$, with $J=J(M)=\min\left\{j\in\N: \; |M|\leq\varepsilon_{j}^{-\sigma}\right\}$.
%Assume that $h$ is ${\cal C}^2$ on each connected component of $\Gamma^{(M)}$, satisfying
% \begin{itemize}
%   \item[\rm (c1)] $\left\{\begin{array}{l}
%   |h|_{\Gamma^{(M)}_j}\leq D \ with \ some \ 0<D <2, \; 0\leq j\leq J+1, \\[2mm]
%   for \; \nu= 1,2 ,  \; |\partial^{\nu}h|_{\Gamma^{(M)}_0}\leq \varepsilon_0^\frac16; \;  |\partial^{\nu} h| \leq \varepsilon_j^{\frac\sigma3}|\xi^{(M)}|^{5-2\nu} \ on \  \Gamma^{(M)}_{j+1}, \ 0\leq j\leq J,
%   \end{array}\right.$,
%   \item[\rm (c2)] for any connected component $(E_*, E_{**})\subset \Gamma^{(M)}_{J+1}$, $\left|\left. h\right|^{E^+_{**}}_{E^-_*}\right|\leq \varepsilon_{J}^{\sigma(1+\frac\sigma6)}$.
% \end{itemize}
For $l=0,1$, we have
\begin{eqnarray*}
\left|\int_{[\inf\sigma(H), \infty)\atop {\rho\leq \varepsilon_0^{-\frac{\sigma}4}} } (h_l-\delta_{l,0}) \, \rho^{2l} \,  \cos (M \rho) \, \partial \rho\, dE\right| &\leq&\frac{\varepsilon_0^{\frac{\sigma^2}{15}}}{|M|^{1+\frac{\sigma}{15}}},  \\
\left|\int_{[\inf\sigma(H), \infty)\atop {\rho> \varepsilon_0^{-\frac{\sigma}4}} } \frac{(h_l-\delta_{l,0}) \,  \rho^{2l}  \cos(M \rho)}{1+\rho^8} \, \partial \rho\, dE\right|&\leq&\frac{\varepsilon_0^{\frac{\sigma^2}{15}}}{|M|^{1+\frac{\sigma}{15}}},
\end{eqnarray*}
and for $k=2l, \, 2l+2$,
\begin{eqnarray*}
\left|\int_{[\inf\sigma(H), \infty)\atop {\rho\leq \varepsilon_0^{-\frac{\sigma}4}} } h_l \, \rho^{k} \,  \cos (M \rho) \, \partial \rho\, dE\right|&\leq&  \frac{\varepsilon_0^{\frac{\sigma^2}{15}}+ 2 \varepsilon_0^{(2l-\frac{k}4)\sigma}}{|M|^{1+\frac{\sigma}{15}}},\\
\left|\int_{[\inf\sigma(H), \infty)\atop {\rho> \varepsilon_0^{-\frac{\sigma}4}} } \frac{h_l \,  \rho^k  \cos(M \rho)}{1+\rho^8} \, \partial \rho\, dE\right|&\leq&\frac{\varepsilon_0^{\frac{\sigma^2}{15}}}{|M|^{1+\frac{\sigma}{15}}}.
\end{eqnarray*}
\end{Proposition}

\begin{Remark}
The initial attempt was to bound the above integral by $\sim\frac{1}{M^2}$, which guarantees the integrability with respect to $M$ over $\{M\in\R:|M|>1\}$. However, we could not get this by a simple application of integration by parts since here $h$ is just piecewise ${\cal C}^2$ on $\Gamma^{(M)}$ and it is not continuous at the edge points.
For this reason, we expect the bound $\sim\frac{1}{|M|^{1+\frac\sigma{15}}}$ instead.
It also guarantees the convergence of the integral on an unbounded interval, which will be applied in the next section.

%On each connected component $(E_*, E_{**})\subset\Gamma_{J+1}^{(M)}$, where $\sin M\rho$ does not vanish at both edge points,
%$\left.h_{l} \right|_{(E_*, E_{**})}$ is well estimated by the ${\cal C}^2$ property.
%As for the external variation of $h$, i.e., to control $\left.h_{l} \right|_{E^-_*}^{E^+_{**}}$(which is also necessary in the integration by parts, as shown in (\ref{ineq_int11})) and (\ref{ineq_int12})), we need an additional condition {\rm (c2)}.
%This is related to (\ref{edge_point}) in {\bf (S4)} of Proposition \ref{propHakan1} and the last statement of Lemma \ref{coefficient_abc}.
\end{Remark}

\noindent{\bf Proof of Proposition \ref{hmn_integral}:}
We consider the above integrals in a more general sense, i.e., to estimate the integrals
$$\int_{[\inf\sigma(H), \infty)\atop {\rho\leq \varepsilon_0^{-\frac{\sigma}4}} } f \cdot \rho^k \cos(M \rho)\, \partial \rho\, dE, \quad \int_{[\inf\sigma(H), \infty)\atop {\rho> \varepsilon_0^{-\frac{\sigma}4}} } \frac{f\cdot \rho^k \cos(M \rho)}{1+\rho^8}\, \partial \rho\, dE, $$
with $f$ ${\cal C}^2$ on each connected component $(E_*, E_{**})\subset\Gamma^{(M)}$, satisfying
\begin{equation}\label{esti_f}
|f|_{\Gamma^{(M)}}\leq 2; \quad |\partial^\nu f|_{\Gamma_0^{(M)}}\leq \varepsilon_0^{\frac16}, \;\ |\partial^\nu f|_{\Gamma_{j+1}^{(M)}}\leq \varepsilon_j^{\frac\sigma3}|\xi^{(M)}|^{5-2\nu},\;\ \nu=1,2; \quad \left|\left. f\right|_{E_*^-}^{E_{**}^+}\right|\leq \varepsilon_J^{\sigma(1+\frac\sigma6)}.
\end{equation}

The above integrals are sums of integrals over the connected component $(E_*,\, E_{**})\subset \Gamma^{(M)}$.
By integration by parts, we have, for $k=0,2,4$,
\begin{eqnarray*}
&&\int_{[\inf\sigma(H), \infty)\atop {\rho\leq \varepsilon_0^{-\frac{\sigma}4}} } f \cdot \rho^k \cos(M \rho)\, \partial \rho\, dE\\
&=&\frac{1}{M}\sum_{(E_*, \, E_{**})\subset \Gamma^{(M)}\atop {\rho\leq \varepsilon_0^{-\frac{\sigma}4}} }\left. f\cdot \rho^k \sin M \rho  \right|_{(E_*, E_{**})} - \frac{1}{M}\sum_{(E_*, \, E_{**})\subset \Gamma^{(M)}\atop {\rho\leq \varepsilon_0^{-\frac{\sigma}4}}}\int_{E_*}^{E_{**} } \partial f\cdot \rho^k\sin M \rho \,  dE\\
& &-\,\frac{1}{M}\sum_{(E_*, \, E_{**})\subset \Gamma^{(M)}\atop {\rho\leq \varepsilon_0^{-\frac{\sigma}4}}}\int_{E_*}^{E_{**}} f\cdot\partial_\rho (\rho^k) \sin M \rho\, \partial \rho \,  dE\\
&=:&{\cal I}_{k,1}+{\cal I}_{k,2}+{\cal I}_{k,3},\\[2mm]
&&\int_{[\inf\sigma(H), \infty)\atop {\rho> \varepsilon_0^{-\frac{\sigma}4}} } \frac{f\cdot \rho^k \cos(M \rho)}{1+\rho^8}\, \partial \rho\, dE\\
&=& \frac{1}{M}\sum_{(E_*, \, E_{**})\subset \Gamma^{(M)}\atop {\rho> \varepsilon_0^{-\frac{\sigma}4}} }\left. \frac{f\cdot \rho^k \sin M \rho }{1+\rho^8}  \right|_{(E_*, E_{**})}- \frac{1}{M}\sum_{(E_*, \, E_{**})\subset \Gamma^{(M)}\atop {\rho> \varepsilon_0^{-\frac{\sigma}4}}}\int_{E_*}^{E_{**} } \frac{\partial f\cdot \rho^k\sin M \rho}{1+\rho^8} \,  dE,\\
& & - \, \frac{1}{M}\sum_{(E_*, \, E_{**})\subset \Gamma^{(M)}\atop {\rho> \varepsilon_0^{-\frac{\sigma}4}}}\int_{E_*}^{E_{**}} f\cdot\partial_\rho \left(\frac{\rho^k}{1+\rho^8}\right) \sin M \rho\, \partial \rho \,  dE \\
&=:&{\cal J}_{k,1}+{\cal J}_{k,2}+{\cal J}_{k,3}.
\end{eqnarray*}
Obviously, ${\cal I}_{0,3}=0$. In this proof, ``$\displaystyle \sum_{(E_*, \, E_{**})\subset \Gamma^{(M)}}$" denotes the sum over all the connected components of $\Gamma^{(M)}$ and it is similar for ``$\displaystyle \sum_{(E_*, \, E_{**})\subset \Gamma_j^{(M)}}$", $0\leq j \leq J+1$. Note that for the last connected component $(E_*, E_{**})\subset\Gamma^{(M)}_{0}$, we have $E_{**}=\infty$.

The proof of proposition comes with the following three lemmas.

\begin{Lemma}\label{lem001} For $k=0,2,4$,
$|{\cal I}_{k,1}|\leq \frac{ \varepsilon_{J}^{\frac{\sigma^2}{12}}+ \varepsilon_{J}^{\frac{\sigma^2}{5}}\varepsilon_0^{-\frac{k\sigma}4}|f|_{\Gamma^{(M)}}}{|M|^{1+\frac\sigma{15}}}$,  $|{\cal J}_{k,1}| \leq \frac{ \varepsilon_{J}^{\frac{\sigma^2}{14}}}{|M|^{1+\frac\sigma{15}}}$.
\end{Lemma}
\proof
Since {\bf (S4)} of Proposition \ref{propHakan1} implies that $\sin M\rho(E)=0$ if $E\in\partial\Gamma^{(M)}\setminus\partial\Gamma_{J+1}^{(M)}$, we can see
\begin{eqnarray}
 && \sum_{(E_*, E_{**})\subset \Gamma^{(M)}\atop {\rho\leq \varepsilon_0^{-\frac{\sigma}4}} }\left. f\cdot\rho^k \sin M \rho\right|_{(E_*, E_{**})} \nonumber\\
&=& \sum_{(E_*, E_{**})\subset\Gamma^{(M)}_{J+1}\atop {\rho\leq \varepsilon_0^{-\frac{\sigma}4}} } \left(
\left. f\cdot\rho^k \sin M \rho\right|_{(E_*, E_{**})} -  \left. f\cdot\rho^k \sin M \rho \right|_{E^-_*}^{E^+_{**}}\right), \label{sum_of_integral_I1}\\
 && \sum_{(E_*, E_{**})\subset \Gamma^{(M)}\atop {\rho> \varepsilon_0^{-\frac{\sigma}4}}}\left. \frac{f\cdot\rho^k \sin M \rho }{1+\rho^8}\right|_{(E_*, E_{**})} \nonumber\\
&=& \sum_{(E_*, E_{**})\subset\Gamma^{(M)}_{J+1}\atop {\rho> \varepsilon_0^{-\frac{\sigma}4}}} \left(
\left. \frac{f\cdot\rho^k \sin M \rho }{1+\rho^8}\right|_{(E_*, E_{**})} -  \left. \frac{f\cdot\rho^k \sin M \rho }{1+\rho^8}\right|_{E^-_*}^{E^+_{**}}\right). \label{sum_of_integral_J1}
\end{eqnarray}
By {\bf (S4)} we know $\left|\left. \rho\right|_{(E_*, E_{**})}\right|\leq 2 \varepsilon_{J}^{\sigma(1+\frac\sigma2)}$,
so, for $k=0,2,4$,
\begin{eqnarray*}
\left|\left. \rho^k \sin M \rho\right|_{(E_*, E_{**})}\right|&\leq& 4|M| \varepsilon_0^{-\frac{k\sigma}4}\varepsilon_{J}^{\sigma(1+\frac\sigma2)} \; {\rm for}  \ \rho\leq \varepsilon_0^{-\frac{\sigma}4}, \\
\left|\left.\frac{\rho^k \sin M \rho}{1+\rho^8}\right|_{(E_*, E_{**})}\right|&\leq& 4|M|\varepsilon_0^{\sigma}\varepsilon_{J}^{\sigma(1+\frac\sigma2)}  \; {\rm for}  \ \rho> \varepsilon_0^{-\frac{\sigma}4}.
\end{eqnarray*}
Indeed, we have
\begin{eqnarray*}
\partial_\rho\left(\rho^k \sin M \rho\right)&=&\left\{\begin{array}{ll}
\displaystyle M \cos(M \rho), & k=0 \\[3mm]
 \displaystyle  k \rho^{k-1}\sin M \rho+ \rho^k  M \cos(M \rho),                                                                                               & k=2,4
\end{array}
 \right. , \\
 \partial_\rho\left(\frac{\rho^k \sin M \rho}{1+\rho^8}\right)&=&\left\{\begin{array}{ll}
\displaystyle \frac{ M \cos(M \rho)}{1+\rho^8}
-\frac{8\rho^{7} \sin M \rho}{(1+\rho^8)^2}, & k=0 \\[3mm]
 \displaystyle  \frac{k \rho^{k-1}\sin M \rho}{1+\rho^8}
+\frac{\rho^k  M \cos(M \rho)}{1+\rho^8}
-\frac{8\rho^{k+7} \sin M \rho}{(1+\rho^8)^2},& k=2,4
\end{array}
 \right. ,
\end{eqnarray*}
which implies $\left|\partial_\rho\left(\rho^k \sin M \rho\right)\right|\leq 2|M| \varepsilon_0^{-\frac{k\sigma}4}$ for $\rho\leq \varepsilon_0^{-\frac{\sigma}4}$, and $\left|\partial_\rho\left(\frac{\rho^k \sin M \rho}{1+\rho^8}\right)\right|\leq 2|M|\varepsilon_0^{\sigma}$  for $\rho> \varepsilon_0^{-\frac{\sigma}4}$. Hence, for $k=0,2,4$,
\begin{eqnarray}
\left|\left. f\cdot \rho^k \, \sin M \rho \right|^{E^+_{**}}_{E^-_*}\right| &\leq&
\left|\left. f \right|^{E^+_{**}}_{E^-_*}\right| \left| \rho^k(E_{**}) \, \sin M \rho(E_{**})\right|+|f(E_{*}^-)|\left|\left.  \rho^k \, \sin M \rho\right|_{(E_*, E_{**})}\right|\nonumber\\
&\leq& \varepsilon_{J}^{\sigma(1+\frac\sigma6)}\cdot \varepsilon_0^{-\frac{k\sigma}4} + |f|_{\Gamma^{(M)}}\cdot 2|M| \varepsilon_0^{-\frac{k\sigma}4}\varepsilon_{J}^{\sigma(1+\frac\sigma2)}\nonumber\\
&\leq&\frac{\varepsilon_J^{\frac{\sigma^2}{10}}}{|M|^{\frac{\sigma}{15}}}+\frac{\varepsilon_J^{\frac{\sigma^2}{4}}\varepsilon_0^{-\frac{k\sigma}4}|f|_{\Gamma^{(M)}}}{|M|^{\frac{\sigma}{5}}} \; {\rm for}  \ \rho\leq \varepsilon_0^{-\frac{\sigma}4},\label{ineq_int11}\\
\left|\left. \frac{f\cdot \rho^k \, \sin M \rho }{1+\rho^8}\right|^{E^+_{**}}_{E^-_*}\right| &\leq&
\left|\left. f \right|^{E^+_{**}}_{E^-_*}\right| \left| \frac{\rho^k \, \sin M \rho }{1+\rho^8}(E_{**})\right|+|f(E_{*}^-)|\left|\left. \frac{ \rho^k \, \sin M \rho}{1+\rho^8}\right|_{(E_*, E_{**})}\right|\nonumber\\
&\leq& \varepsilon_{J}^{\sigma(1+\frac\sigma6)}\cdot \varepsilon_0^{\sigma} + |f|_{\Gamma^{(M)}}\cdot 2|M|\varepsilon_0^{\sigma}\varepsilon_{J}^{\sigma(1+\frac\sigma2)}\nonumber\\
&\leq&\frac{\varepsilon_J^{\frac{\sigma^2}{10}}}{|M|^{\frac{\sigma}{15}}}+\frac{\varepsilon_J^{\frac{\sigma^2}{4}}\varepsilon_0^{\sigma}}{|M|^{\frac{\sigma}{5}}} \; {\rm for}  \ \rho> \varepsilon_0^{-\frac{\sigma}4} ,\label{ineq_int12}
\end{eqnarray}
and similarly, since $E_{**}-E_{*}\leq \varepsilon_{J}^{\sigma(1+\frac{\sigma}{3})}$ for the connected component $(E_{*},E_{**})\subset \Gamma_{J+1}^{(M)}$,
\begin{eqnarray*}
\left|\left. f\cdot \rho^k \, \sin M \rho \right|_{(E_*,E_{**})}\right|
&\leq&  \varepsilon_J^{\sigma(1+\frac{\sigma}3)} \cdot \varepsilon_0^{-\frac{k\sigma}4} + |f|_{\Gamma^{(M)}} \cdot 2|M| \varepsilon_0^{-\frac{k\sigma}4}\varepsilon_{J}^{\sigma(1+\frac\sigma2)}\\
&\leq& \frac{\varepsilon_J^{\frac{\sigma^2}{10}}}{|M|^{\frac{\sigma}{15}}}+\frac{\varepsilon_J^{\frac{\sigma^2}{4}}\varepsilon_0^{-\frac{k\sigma}4}|f|_{\Gamma^{(M)}}}{|M|^{\frac{\sigma}{5}}} \; {\rm for}  \ \rho\leq \varepsilon_0^{-\frac{\sigma}4} ,\\
\left|\left. \frac{f\cdot \rho^k \, \sin M \rho }{1+\rho^8}\right|_{(E_*,E_{**})}\right|
%&\leq&
%\left|\left. f \right|_{(E_*,E_{**})}\right| \left| \frac{\rho^k \, \sin M \rho }{1+\rho^8}(E_{**})\right|+|f(E_{*}^+)|\left|\left. \frac{ \rho^k \, \sin M \rho}{1+\rho^8}\right|_{(E_*, E_{**})}\right|\\
&\leq& \varepsilon_{J}^{\sigma(1+\frac\sigma3)}\cdot \varepsilon_0^{\sigma} +|f|_{\Gamma^{(M)}}\cdot 2|M|\varepsilon_0^{\sigma}\varepsilon_{J}^{\sigma(1+\frac\sigma2)} \\
&\leq&\frac{\varepsilon_J^{\frac{\sigma^2}{10}}}{|M|^{\frac{\sigma}{15}}}+\frac{\varepsilon_J^{\frac{\sigma^2}{4}}\varepsilon_0^{\sigma}}{|M|^{\frac{\sigma}{5}}} \; {\rm for}  \ \rho> \varepsilon_0^{-\frac{\sigma}4}.
\end{eqnarray*}
Recalling that there are at most $|\ln\varepsilon_0|^{(J+1)^3 d}$ connected components in $[\inf\sigma(H), \infty)$, so
\begin{eqnarray*}
\left|\frac{1}{M}\sum_{(E_*, E_{**})\subset\Gamma^{(M)}\atop {\rho\leq \varepsilon_0^{-\frac{\sigma}4}}}\left. h\cdot\rho^k \sin M \rho\right|_{(E_*, E_{**})}\right|&\leq& \frac{|\ln\varepsilon_0|^{(J+1)^3d}}{|M|} \cdot \frac{2\varepsilon_J^{\frac{\sigma^2}{10}}+2\varepsilon_J^{\frac{\sigma^2}{4}}\varepsilon_0^{-\frac{k\sigma}4}|f|_{\Gamma^{(M)}}}{|M|^{\frac{\sigma}{15}}}\\
&\leq&\frac{ \varepsilon_{J}^{\frac{\sigma^2}{12}}+ \varepsilon_{J}^{\frac{\sigma^2}{5}}\varepsilon_0^{-\frac{k\sigma}4}|f|_{\Gamma^{(M)}}}{|M|^{1+\frac\sigma{15}}},
\end{eqnarray*}
$$\left|\frac{1}{M}\sum_{(E_*, E_{**})\subset\Gamma^{(M)}\atop {\rho> \varepsilon_0^{-\frac{\sigma}4}}}\left. \frac{h\cdot\rho^k \sin M \rho }{1+\rho^8}\right|_{(E_*, E_{**})}\right|\leq \frac{|\ln\varepsilon_0|^{(J+1)^3d}}{|M|} \cdot \frac{2\varepsilon_J^{\frac{\sigma^2}{10}}+2\varepsilon_J^{\frac{\sigma^2}{4}}\varepsilon_0^{\sigma}}{|M|^{\frac{\sigma}{15}}}
\leq\frac{\varepsilon_{J}^{\frac{\sigma^2}{14}}}{|M|^{1+\frac\sigma{15}}}.$$
\qed

\begin{Lemma}\label{lem002}
$|{\cal I}_{k,2}|, \, |{\cal J}_{k,2}| \leq \frac{\varepsilon_0^{\frac{10\sigma}9}}{4|M|^{\frac32}}$ for $k=0,2,4$.
\end{Lemma}
\proof Let $\rho^{(M)}:=\Re\alpha^{(M)}$. By {\bf (S1)} in Proposition \ref{propHakan1}, we have $|\rho-\rho^{(M)}|_{\Gamma^{(M)}}<\varepsilon_J^{\frac14}$.
On each connected component $(E_*, E_{**})\subset\Gamma^{(M)}$, to estimate ${\cal I}_{k,2}$ and ${\cal J}_{k,2}$, we can consider the approximated integral $\displaystyle \int_{E_*}^{E_{**}} \, \partial f\cdot(\rho^{(M)})^k\sin M \rho^{(M)} \,  dE$ and
$\displaystyle \int_{E_*}^{E_{**}} \, \frac{\partial f\cdot(\rho^{(M)})^k\sin M \rho^{(M)}}{1+(\rho^{(M)})^8}\,  dE$ instead. Indeed, by noting that $|M|\leq \varepsilon_J^{-\sigma}$, it is easy to verify that, for $k=0,2,4$,
\begin{eqnarray*}
\left|(\rho^{(M)})^k\, \sin  M \rho^{(M)}-\rho^k\, \sin M \rho\right|&\leq& \varepsilon_J^{\frac15} \;\ {\rm for} \ \rho\leq\varepsilon_0^{-\frac\sigma4}, \\
\left|\frac{(\rho^{(M)})^k\, \sin  M \rho^{(M)}}{1+(\rho^{(M)})^8}-\frac{\rho^k\, \sin M \rho}{1+\rho^8}\right|&\leq& \frac{\varepsilon_J^{\frac15}}{1+\rho^8} \;\ {\rm for} \ \rho>\varepsilon_0^{-\frac\sigma4}.
\end{eqnarray*}
Then, by the transversality of $\rho$ in (\ref{transversality_rho}), the gap estimate in (\ref{gap-estime}), and the fact that $|\partial f|_{\Gamma^{(M)}}\leq \varepsilon_0^{2\sigma}$ in view of (\ref{esti_f}), we have
\begin{eqnarray*}
 & & \sum_{(E_*, E_{**})\subset\Gamma^{(M)}\atop {\rho\leq \varepsilon_0^{-\frac{\sigma}4}}}\int_{E_*}^{E_{**}} \, | \partial f|\left|(\rho^{(M)})^k\, \sin  M \rho^{(M)}-\rho^k\, \sin M \rho\right| \,  dE \\
&\leq& \varepsilon_0^{2\sigma}\varepsilon_J^{\frac15}  \int_{[\inf\sigma(H), \infty)\atop {\rho\leq \varepsilon_0^{-\frac{\sigma}4}}}  \,  dE\\
&\leq&2 \varepsilon_0^{2\sigma}\varepsilon_J^{\frac15}\int_{0}
^{\varepsilon_0^{-\frac{\sigma}4}} \rho \,  d\rho+ c' \varepsilon_0^{2\sigma}\varepsilon_J^{\frac15}\sum_{k\in\Z^d} e^{-\beta'|k|^{\iota}}\\
&\leq& \varepsilon_0^{\sigma}\varepsilon_J^{\frac15} \;\ {\rm for} \ \rho\leq\varepsilon_0^{-\frac\sigma4},\\[1mm]
 & & \sum_{(E_*, E_{**})\subset\Gamma^{(M)}\atop {\rho> \varepsilon_0^{-\frac{\sigma}4}}}\int_{E_*}^{E_{**}} \, | \partial f| \left|\frac{(\rho^{(M)})^k\, \sin  M \rho^{(M)}}{1+(\rho^{(M)})^8}-\frac{\rho^k\, \sin M \rho}{1+\rho^8}\right| \,  dE \\
&\leq& \varepsilon_0^{2\sigma}\varepsilon_J^{\frac15}\int_{[\inf\sigma(H),\infty) \atop {\rho> \varepsilon_0^{-\frac{\sigma}4}} } \frac{1}{1+\rho^8} \,  dE\\
&\leq&2 \varepsilon_0^{2\sigma}\varepsilon_J^{\frac15}\int_{\varepsilon_0^{-\frac{\sigma}4}}^{\infty} \frac{\rho}{1+\rho^8} \,  d\rho
+c' \varepsilon_0^{2\sigma}\varepsilon_J^{\frac15}\sum_{k\in\Z^d} e^{-\beta'|k|^{\iota}}\\
&\leq& \varepsilon_0^{\sigma}\varepsilon_J^{\frac15} \;\ {\rm for} \ \rho>\varepsilon_0^{-\frac\sigma4}.
\end{eqnarray*}
Note that we consider the above integrals over $[\inf\sigma(H),\infty)$ in two parts.
\begin{itemize}
\item In $\sigma(H)$, we have the transversality (\ref{transversality_rho}) of $\rho$, so this is transformed into an integral with respect to $\rho$.
\item In the gaps of spectrum, we have $\partial\rho=0$, so the function of integral is a constant. This part is controlled by the gap estimate in (\ref{gap-estime}).
\end{itemize}

Now we consider the approximated integrals $\displaystyle \int_{E_*}^{E_{**}} \, \partial f\cdot(\rho^{(M)})^k\sin M \rho^{(M)} \,  dE$
on the connected component $(E_*, E_{**})\subset\Gamma^{(M)}$ where $\rho\leq \varepsilon_0^{-\frac{\sigma}4}$,
and $\displaystyle \int_{E_*}^{E_{**}} \, \frac{\partial f\cdot(\rho^{(M)})^k\sin M \rho^{(M)}}{1+(\rho^{(M)})^8}\,  dE$
on the connected component where $\rho> \varepsilon_0^{-\frac{\sigma}4}$,
\begin{itemize}
  \item On $(E_*, E_{**})\subset\Gamma_0^{(M)}$, we have $\rho^{(M)}=\xi^{(M)}$.
So, to compute the integrals $\displaystyle \int_{E_*}^{E_{**}} \partial f\cdot(\rho^{(M)})^k \sin M \rho^{(M)} \,  dE$ and
$\displaystyle \int_{E_*}^{E_{**}} \frac{\partial f\cdot(\rho^{(M)})^k \sin M \rho^{(M)}}{1+(\rho^{(M)})^8} \,  dE$, we assume $\xi^{(M)}\neq0$.
Hence, by {\bf (S2)}, we have
  $\partial\rho^{(M)}=\frac{\partial\det A^{(M)}}{2\rho^{(M)}}> \frac1{3\rho^{(M)}}$. Then
\begin{eqnarray*}
\int_{E_*}^{E_{**}} \partial f\cdot(\rho^{(M)})^k \sin M \rho^{(M)} \,  dE
&=& 2\int_{E_*}^{E_{**}} \frac{\partial f\cdot(\rho^{(M)})^{k+1}}{\partial\det A^{(M)}}\sin  M \rho^{(M)}\, \partial \rho^{(M)}  dE\\
&=& -\frac{2}{M}\left[\left.\frac{\partial f\cdot(\rho^{(M)})^{k+1} \cos  M \rho^{(M)}}{\partial\det A^{(M)}}\right|_{(E_*, E_{**})}\right.\\
& & - \, \left. \int_{E_*}^{E_{**}} \partial\left(\frac{\partial f\cdot(\rho^{(M)})^{k+1}}{\partial\det A^{(M)}}\right) \cos  M \rho^{(M)} dE\right],
\end{eqnarray*}
\begin{eqnarray*}
&& \int_{E_*}^{E_{**}} \frac{\partial f\cdot(\rho^{(M)})^k \sin M \rho^{(M)}}{1+(\rho^{(M)})^8} \,  dE \\
&=& 2\int_{E_*}^{E_{**}} \frac{\partial f\cdot(\rho^{(M)})^{k+1}\sin  M \rho^{(M)}}{\partial\det A^{(M)}\cdot(1+(\rho^{(M)})^8)}\partial \rho^{(M)}  dE\\
&=& -\frac{2}{M}\left.\frac{\partial f\cdot(\rho^{(M)})^{k+1} \cos  M \rho^{(M)}}{\partial\det A^{(M)}\cdot(1+(\rho^{(M)})^8)}\right|_{(E_*, E_{**})}\\
& & + \, \frac{2}{M} \int_{E_*}^{E_{**}} \partial\left(\frac{\partial f\cdot(\rho^{(M)})^{k+1}}{\partial\det A^{(M)}\cdot(1+(\rho^{(M)})^8)}\right) \cos  M \rho^{(M)} dE.
\end{eqnarray*}
(\ref{sigma_m_0}) implies that $|\partial\det A^{(M)}-1|$, $|\partial^2\det A^{(M)}| \leq 2 \varepsilon_0^{\frac23}$, then for $\rho^{(M)}\leq\varepsilon_0^{-\frac\sigma4}+\varepsilon_J^{\frac14}$,
$\displaystyle \left|\frac{\partial f\cdot(\rho^{(M)})^{k+1}}{\partial\det A^{(M)}} \right|\leq \varepsilon_0^{\frac17-\frac{5\sigma}{4}}$,
\begin{eqnarray*}
 \left|\partial\left(\frac{\partial f\cdot(\rho^{(M)})^{k+1}}{\partial\det A^{(M)}}\right)\right|
&\leq&|\rho^{(M)}|^{k+1}\cdot\left[\frac{|\partial^2 f|}{\partial\det A^{(M)}}+\frac{|\partial f|\cdot|\partial^2\det A^{(M)}|}{(\partial\det A^{(M)})^2}\right]\\
&& + \,\frac{|\partial f|\cdot(k+1)(\rho^{(M)})^{k}\cdot \partial\rho^{(M)}}{\partial\det A^{(M)}} \\
&\leq&2\varepsilon_0^{\frac16} \cdot(3(\rho^{(M)})^{k+2}+(k+1)(\rho^{(M)})^{k})\partial \rho^{(M)}\\
&\leq& \varepsilon_0^{\frac17} \cdot(\rho^{(M)})^{k+2} \partial \rho^{(M)},
\end{eqnarray*}
and for $\rho^{(M)}>\varepsilon_0^{-\frac\sigma4}-\varepsilon_J^{\frac14}$, $\displaystyle  \left|\frac{\partial f\cdot(\rho^{(M)})^{k+1}}{\partial\det A^{(M)}\cdot(1+(\rho^{(M)})^8)} \right|\leq \varepsilon_0^{\frac17+\sigma}$,
\begin{eqnarray*}
&&\left|\partial\left(\frac{\partial f\cdot(\rho^{(M)})^{k+1}}{\partial\det A^{(M)}\cdot(1+(\rho^{(M)})^8)}\right)\right|\\
&\leq& \frac{|\rho^{(M)}|^{k+1}}{1+(\rho^{(M)})^8}\cdot\left[\frac{|\partial^2 f|}{\partial\det A^{(M)}}+\frac{|\partial f|\cdot|\partial^2\det A^{(M)}|}{(\partial\det A^{(M)})^2}\right]\\
&&+\frac{|\partial f|\cdot \partial\rho^{(M)}}{\partial\det A^{(M)}}\cdot \left[\frac{(k+1)(\rho^{(M)})^{k}}{1+(\rho^{(M)})^8} + \frac{8(\rho^{(M)})^{k+8}}{(1+(\rho^{(M)})^8)^2}\right]  \\
&\leq&\varepsilon_0^{\frac17} \frac{(\rho^{(M)})^{14}\, \partial \rho^{(M)}}{(1+(\rho^{(M)})^8)^2}
\end{eqnarray*}
Therefore, for $(E_{*}, E_{**})\subset \Gamma^{(M)}_0$,
\begin{eqnarray*}
\left|\int_{E_*}^{E_{**}} \partial f\cdot (\rho^{(M)})^k  \sin M \rho^{(M)} \,  dE\right|&\leq&\frac{\varepsilon_0^{\frac18}}{|M|}\;\ {\rm for} \ \rho^{(M)}\leq\varepsilon_0^{-\frac\sigma4}+\varepsilon_J^{\frac14},  \\
\left|\int_{E_*}^{E_{**}} \frac{\partial f\cdot (\rho^{(M)})^k  \sin M \rho^{(M)}}{1+(\rho^{(M)})^8} \,  dE\right|&\leq&\frac{\varepsilon_0^{\frac18}}{|M|}\;\ {\rm for} \ \rho^{(M)}>\varepsilon_0^{-\frac\sigma4}-\varepsilon_J^{\frac14},
\end{eqnarray*}
by combining the above estimates, and noting that on the connected component $(E_*, E_{**})$ where $-\varepsilon_J^{\frac14}\leq \rho^{(M)}\leq\varepsilon_0^{-\frac\sigma4}+\varepsilon_J^{\frac14}$,
$$\int_{E_*}^{E_{**}} \left|\partial\left(\frac{\partial f\cdot(\rho^{(M)})^{k+1}}{\partial\det A^{(M)}}\right)\right| dE
   \leq \varepsilon_0^{\frac17} \int_{-\varepsilon_J^{\frac14}}^{\varepsilon_0^{-\frac\sigma4}+\varepsilon_J^{\frac14}} (\rho^{(M)})^{k+2} \, d \rho^{(M)}
\leq \frac14\varepsilon_0^{\frac{1}{8}},$$
and on the connected component $(E_*, E_{**})$ where $\rho^{(M)}>\varepsilon_0^{-\frac\sigma4}-\varepsilon_J^{\frac14}$,
\begin{eqnarray*}
\int_{E_*}^{E_{**}} \left|\partial\left(\frac{\partial f\cdot(\rho^{(M)})^{k+1}}{\partial\det A^{(M)}(1+(\rho^{(M)})^8)}\right)\right| dE
&\leq&  \varepsilon_0^{\frac17} \int_{\varepsilon_0^{-\frac{\sigma}4}-\varepsilon_J^{\frac14}}^\infty  \frac{(\rho^{(M)})^{14}\, \partial \rho^{(M)}}{(1+(\rho^{(M)})^8)^2}\partial \rho^{(M)} dE\\
&\leq& \frac14\varepsilon_0^{\frac18}.
\end{eqnarray*}

  \item On $(E_*, \,  E_{**})\subset\Gamma_{j+1}^{(M)}$, recall that there is an interval ${\cal I}\subset (E_*, \,  E_{**})$, such that $\xi^{(M)}=0$.
      So (\ref{esti_f}) implies $\partial f=0$ on ${\cal I}$. On $(E_*, \,  E_{**})\setminus {\cal I}$, noting that $\partial\rho^{(M)}=\partial\xi^{(M)}$ and in view of (\ref{esti_plat}), we have for $-\varepsilon_J^{\frac14} \leq \rho^{(M)}\leq\varepsilon_0^{-\frac\sigma4}+\varepsilon_J^{\frac14}$,
      $$\left| \frac {\partial f\cdot (\rho^{(M)})^k}{\partial \rho^{(M)}} \right|, \ \left|\partial\left(\frac {\partial f\cdot (\rho^{(M)})^k}{ \partial \rho^{(M)}}\right)\right|\leq \frac14\varepsilon_0^{\frac{9\sigma}8},\quad k=0,2,4,$$
and for $\rho^{(M)}>\varepsilon_0^{-\frac\sigma4}-\varepsilon_J^{\frac14}$,
$$\left| \frac {\partial f\cdot (\rho^{(M)})^k}{(1+(\rho^{(M)})^8)\partial \rho^{(M)}} \right|, \;\  \left|\partial\left(\frac {\partial f\cdot (\rho^{(M)})^k}{(1+(\rho^{(M)})^8) \partial \rho^{(M)}}\right)\right| \leq \varepsilon_0^{2\sigma},\quad k=0,2,4 .$$
Then, with ${\cal P}_1$ and ${\cal P}_2$ denoting the two connected components of $(E_*, \,  E_{**})\setminus {\cal I}$,
\begin{eqnarray*}
&&\int_{(E_*, \,  E_{**})\setminus {\cal I}} \, \partial f\cdot (\rho^{(M)})^k\sin  M \rho^{(M)} \,  dE\\
&=&\int_{(E_*, \,  E_{**})\setminus {\cal I}} \, \frac{\partial f\cdot (\rho^{(M)})^k \sin  M \rho^{(M)}}{\partial \rho^{(M)}} \partial \rho^{(M)}\, dE  \\
&=&\frac{-1}{M}\left.\frac{\partial f\cdot (\rho^{(M)})^k}{\partial \rho^{(M)}}\, \cos M\rho^{(M)}\right|_{{\cal P}_1}+\, \frac{-1}{M}\left.\frac{\partial f\cdot (\rho^{(M)})^k}{\partial \rho^{(M)}}\, \cos M\rho^{(M)}\right|_{{\cal P}_2}\\[1mm]
& &+\, \frac{1}{M}\int_{{\cal P}_1\cup{\cal P}_2} \, \partial\left(\frac {\partial f\cdot (\rho^{(M)})^k}{\partial \rho^{(M)}}\right)\, \cos M\rho^{(M)} \,  dE,
\end{eqnarray*}
\begin{eqnarray*}
&&\int_{(E_*, \,  E_{**})\setminus {\cal I}} \, \frac{\partial f\cdot (\rho^{(M)})^k\sin  M \rho^{(M)}}{1+(\rho^{(M)})^8} \,  dE\\
&=&\int_{(E_*, \,  E_{**})\setminus {\cal I}} \, \frac {\partial f\cdot (\rho^{(M)})^k}{(1+(\rho^{(M)})^8)\partial \rho^{(M)}}\sin  M \rho^{(M)} \cdot \partial \rho^{(M)}\, dE  \\
&=&\frac{-1}{M}\left.\frac{\partial h\cdot (\rho^{(M)})^k}{(1+(\rho^{(M)})^8)\partial \rho^{(M)}}\, \cos M\rho^{(M)}\right|_{{\cal P}_1}+ \frac{-1}{M}\left.\frac{\partial h\cdot (\rho^{(M)})^k}{(1+(\rho^{(M)})^8)\partial \rho^{(M)}}\, \cos M\rho^{(M)}\right|_{{\cal P}_2}\\[1mm]
& &+\, \frac{1}{M}\int_{{\cal P}_1\cup{\cal P}_2} \, \partial\left(\frac {\partial h\cdot (\rho^{(M)})^k}{(1+(\rho^{(M)})^8)\partial \rho^{(M)}}\right)\, \cos M\rho^{(M)} \,  dE.
\end{eqnarray*}
By noting that every connected component of $\Gamma^{(M)}_{j+1}$, $j\geq 0$, is shorter than $\varepsilon_0^{\frac{\sigma}{2}}$, we have that both integrals above can be bounded by $\frac{\varepsilon_0^{\frac{9\sigma}8} }{2|M|}$, provided that $\rho^{(M)}\leq\varepsilon_0^{-\frac\sigma4}+\varepsilon_J^{\frac14}$ and $\rho^{(M)}>\varepsilon_0^{-\frac\sigma4}-\varepsilon_J^{\frac14}$ respectively.
\end{itemize}
So, for each $(E_*, E_{**})\subset \Gamma^{(M)}$, we have
\begin{eqnarray*}
&&\left| \int_{E_*}^{E_{**}} \, \partial f\cdot(\rho^{(M)})^k\sin  M \rho^{(M)}\,  dE\right|\leq \frac{\varepsilon_0^{\frac{9\sigma}8} }{2|M|} +\varepsilon_J^{\frac15}\leq  \frac{ \varepsilon_0^{\frac{9\sigma}8}}{|M|} \;\ {\rm for} \ \rho\leq\varepsilon_0^{-\frac\sigma4},  \\[2mm]
&&\left| \int_{E_*}^{E_{**}} \, \frac{\partial f\cdot(\rho^{(M)})^k\sin  M \rho^{(M)}}{1+(\rho^{(M)})^8} \,  dE\right|\leq \frac{\varepsilon_0^{\frac{9\sigma}8} }{2|M|} +\varepsilon_J^{\frac15}\leq  \frac{ \varepsilon_0^{\frac{9\sigma}8}}{|M|} \;\ {\rm for} \ \rho>\varepsilon_0^{-\frac\sigma4},
\end{eqnarray*}
and then
\begin{equation}\label{ele11}
\left| \frac{1}{M} \sum_{(E_*, E_{**})\subset \Gamma^{(M)}\atop {\rho\leq \varepsilon_0^{-\frac{\sigma}4}}}\int_{E_*}^{E_{**}} \, \frac{\partial f\cdot(\rho^{(M)})^k\sin  M \rho^{(M)}}{1+(\rho^{(M)})^8} \,  dE\right|\leq \frac{|\ln\varepsilon_0|^{(J+1)^3 d}\,  \varepsilon_0^{\frac{9\sigma}8} }{M^2}
\leq\frac{\varepsilon_0^{\frac{10\sigma}9}}{4|M|^{\frac32}},
  \end{equation}
  \begin{equation}\label{ele12}
\left| \frac{1}{M} \sum_{(E_*, E_{**})\subset\Gamma^{(M)}\atop {\rho> \varepsilon_0^{-\frac{\sigma}4}}}\int_{E_*}^{E_{**}} \, \frac{\partial f\cdot(\rho^{(M)})^k\sin  M \rho^{(M)}}{1+(\rho^{(M)})^8} \,  dE\right|\leq \frac{\varepsilon_0^{\frac{10\sigma}9}}{4|M|^{\frac32}}.
  \end{equation}
Note that in getting (\ref{ele11}) and (\ref{ele12}), we need to consider two cases about $M$:
\begin{enumerate}
  \item [(\uppercase\expandafter{\romannumeral1})] If $|M|\leq \varepsilon_0^{-\sigma}$, which means $J(M)=0$, then, by noting that $|M|>1$, we get
$$|\ln\varepsilon_0|^{(J+1)^3 d}\, \varepsilon_0^{\frac{9\sigma}8}\leq \frac14 \varepsilon_0^{\frac{10\sigma}9}|M|^{\frac12}.$$
  \item [(\uppercase\expandafter{\romannumeral2})] If $|M|> \varepsilon_0^{-\sigma}$, which means $J(M)\geq 1$ and $|M|> \varepsilon_{J-1}^{-\sigma}=\varepsilon_{0}^{-\sigma(1+\sigma)^{J-1}}$, then
$$\frac{|\ln\varepsilon_0|^{(J+1)^3 d}\, \varepsilon_0^{\frac{9\sigma}8}}{M^2}\leq \frac{|\ln\varepsilon_0|^{(J+1)^3 d}\, \varepsilon_0^{\frac{9\sigma}8}\cdot\varepsilon_{J-1}^{\frac\sigma2}}{|M|^{\frac32}}\leq \frac{\varepsilon_0^{\frac{9\sigma}8}\varepsilon_{J}^{\frac\sigma6}}{4|M|^{\frac32}}.$$\qed
\end{enumerate}

\begin{Lemma}\label{lem003}
$|{\cal I}_{k,3}|\leq \frac{\varepsilon_0^{-\frac{(2k-1)\sigma}8}|f|_{\Gamma^{(M)}}+\varepsilon_0^{\frac{2\sigma}3}}{|M|^\frac32}$ for $k=2,4$, and $|{\cal J}_{k,3}| \leq \frac{\varepsilon_0^{\sigma}}{|M|^{\frac32}}$ for $k=0, 2, 4$.
\end{Lemma}
\proof By the integration by parts, we have
  \begin{eqnarray}
&&-\frac{1}{M}\sum_{(E_*, \, E_{**})\subset \Gamma^{(M)}\atop {\rho\leq \varepsilon_0^{-\frac{\sigma}4}}}\int_{E_*}^{E_{**}} f\cdot\partial_\rho\left(\rho^k\right) \sin M \rho\, \partial\rho \,  dE\nonumber\\
&=& \frac{1}{M^2}\sum_{(E_*,\,  E_{**})\subset \Gamma^{(M)}\atop {\rho\leq \varepsilon_0^{-\frac{\sigma}4}}}\left.f\cdot\partial_\rho\left(\rho^k\right)\cos(M \rho)\right|_{(E_*,\,  E_{**})}\label{part0221-1}\\
  & &  -\frac{1}{M^2}\sum_{(E_*,\,  E_{**})\subset \Gamma^{(M)}\atop {\rho\leq \varepsilon_0^{-\frac{\sigma}4}}}\int_{E_*}^{E_{**}} \partial \left(f\cdot\partial_\rho\left(\rho^k\right) \right) \cos(M \rho)  \, dE,\label{part0222-1}
  \end{eqnarray}
 \begin{eqnarray}
&&-\frac{1}{M}\sum_{(E_*, \, E_{**})\subset \Gamma^{(M)}\atop {\rho> \varepsilon_0^{-\frac{\sigma}4}}}\int_{E_*}^{E_{**}} f\cdot\partial_\rho\left(\frac{\rho^k}{1+\rho^8}\right) \sin M \rho\, \partial\rho \,  dE\nonumber\\
&=& \frac{1}{M^2}\sum_{(E_*,\,  E_{**})\subset \Gamma^{(M)}\atop {\rho> \varepsilon_0^{-\frac{\sigma}4}}}\left.f\cdot\partial_\rho\left(\frac{ \rho^k}{1+\rho^8}\right)\cos(M \rho)\right|_{(E_*,\,  E_{**})}\label{part0221-2}\\
  & &  -\frac{1}{M^2}\sum_{(E_*,\,  E_{**})\subset \Gamma^{(M)}\atop {\rho> \varepsilon_0^{-\frac{\sigma}4}}}\int_{E_*}^{E_{**}} \partial \left(f\cdot\partial_\rho\left(\frac{ \rho^k}{1+\rho^8}\right) \right) \cos(M \rho)  \, dE.\label{part0222-2}
  \end{eqnarray}

In (\ref{part0221-1}) and (\ref{part0221-2}), for each connected component $(E_*,\,  E_{**})\subset\Gamma^{(M)}$, we have
\begin{eqnarray*}
&& \left|\left.f\cdot\partial_\rho\left(\rho^k\right)\cos(M \rho)\right|_{(E_*,\,  E_{**})}\right|\nonumber  \\
&\leq& k|f(E^-_{**})|\left|\left. \rho^{k-1}\cos(M \rho)\right|_{(E_*,\,  E_{**})}\right|
 + k\left|\left.f\right|_{(E_*,\,  E_{**})}\right|\left|\rho^{k-1}(E_*)\cos(M \rho)(E_*)\right| \nonumber  \\
&\leq& 16 |f|_{\Gamma^{(M)}}\varepsilon_0^{-\frac{(k-1)\sigma}4} \;\ {\rm for} \ k=2,4, \  \rho\leq\varepsilon_0^{-\frac{\sigma}{4}},
\end{eqnarray*}
and since
$\left| \partial_\rho\left(\frac{\rho^k}{1+\rho^8}\right) \right| < \varepsilon_0^{\frac{6\sigma}5}$ for $k=0,2,4$, $\rho>\varepsilon_0^{-\frac{\sigma}{4}}$,
\begin{eqnarray*}
 \left|\left.f\cdot\partial_\rho\left(\frac{\rho^k}{1+\rho^8}\right)\cos(M \rho)\right|_{(E_*,\,  E_{**})}\right|\nonumber
&\leq& |f(E^-_{**})|\left|\left. \partial_\rho\left(\frac{\rho^k}{1+\rho^8}\right)\cos(M \rho)\right|_{(E_*,\,  E_{**})}\right|\\
& & + \left|\left.f\right|_{(E_*,\,  E_{**})}\right|\left|\partial_\rho\left(\frac{\rho^k}{1+\rho^8}\right)(E_*)\cos(M \rho)(E_*)\right| \nonumber  \\
&\leq&16 |f|_{\Gamma^{(M)}}\varepsilon_0^{\frac{6\sigma}5} \;\ {\rm for} \ k=0,2,4, \  \rho>\varepsilon_0^{-\frac{\sigma}{4}}.
\end{eqnarray*}
Hence, similar to (\ref{ele11}) and (\ref{ele12}), we get, for $k=2,4$,
\begin{equation}\label{0ele2-1}
\left| \frac{1}{M^2}\sum_{(E_*,\,  E_{**})\subset \Gamma^{(M)}\atop {\rho\leq \varepsilon_0^{-\frac{\sigma}4}}}\left.f\cdot\partial_\rho\left( \rho^k\right)\cos(M \rho)\right|_{(E_*,\,  E_{**})}\right|\leq\frac{\varepsilon_0^{-\frac{(2k-1)\sigma}8}|f|_{\Gamma^{(M)}}}{2|M|^\frac32},
\end{equation}
and for $k=0,2,4$,
\begin{equation}\label{0ele2-2}
\left| \frac{1}{M^2}\sum_{(E_*,\,  E_{**})\subset \Gamma^{(M)}\atop {\rho> \varepsilon_0^{-\frac{\sigma}4}}}\left.f\cdot\partial_\rho\left(\frac{ \rho^k}{1+\rho^8}\right)\cos(M \rho)\right|_{(E_*,\,  E_{**})}\right|\leq\frac{\varepsilon_0^{\sigma}}{2|M|^\frac32}.
\end{equation}

\smallskip

In (\ref{part0222-1}) and (\ref{part0222-2}),
we have
\begin{eqnarray*}
&&\left|-\frac{1}{M^2}\sum_{(E_*,\,  E_{**})\subset \Gamma^{(M)}\atop {\rho\leq \varepsilon_0^{-\frac{\sigma}4}}}\int_{E_*}^{E_{**}} \partial \left(f\cdot\partial_{\rho}\left(\rho^k\right) \right) \cos(M \rho)  \, dE\right|\nonumber\\
&\leq& \frac{1}{M^2}\int_{[\inf\sigma(H),\infty)\atop {\rho\leq \varepsilon_0^{-\frac{\sigma}4}}} |f| \left|\partial_{\rho}^2\left( \rho^k \right)\right| \partial\rho \, dE
+ \frac{1}{M^2}\int_{[\inf\sigma(H),\infty)\atop {\rho\leq \varepsilon_0^{-\frac{\sigma}4}}}|\partial f| \left|\partial_{\rho} \left(\rho^k \right)\right| \, dE,\\[2mm]
&&\left|-\frac{1}{M^2}\sum_{(E_*,\,  E_{**})\subset \Gamma^{(M)}\atop {\rho> \varepsilon_0^{-\frac{\sigma}4}}}\int_{E_*}^{E_{**}} \partial \left(f\cdot\partial_{\rho}\left(\frac{ \rho^k}{1+\rho^8}\right) \right) \cos(M \rho)  \, dE\right|\nonumber\\
&\leq& \frac{1}{M^2}\int_{\inf\sigma(H)}^{\infty} |f| \left|\partial_{\rho}^2\left( \frac{ \rho^k }{1+\rho^8}\right)\right| \partial\rho \, dE
+ \frac{1}{M^2}\int_{\inf\sigma(H)}^{\infty}|\partial f| \left|\partial_{\rho} \left( \frac{\rho^k}{1+\rho^8}\right)\right| \, dE.
\end{eqnarray*}
Since $|\partial f|\leq\varepsilon_0^{2\sigma}$ in each connected component of $\Gamma^{(M)}$, we can see, for $k=2,4$,
$$\frac{1}{M^2}\int_{[\inf\sigma(H),\infty)\atop {\rho\leq \varepsilon_0^{-\frac{\sigma}4}}} |f| \left|\partial_{\rho}^2\left( \rho^k \right)\right| \partial\rho \, dE\leq\frac{\varepsilon_0^{-\frac{(2k-1)\sigma}8}|f|_{\Gamma^{(M)}}}{4M^2},$$
$$ \frac{1}{M^2}\int_{[\inf\sigma(H),\infty)\atop {\rho\leq \varepsilon_0^{-\frac{\sigma}4}}} |\partial f| \left|\partial_{\rho}\left(\rho^k\right)\right|\, dE
\leq\frac{c'\varepsilon_0^{2\sigma}}{M^2} \sum_{l\in\Z^d}e^{-\beta'|l|^\iota}
+\frac{2k\varepsilon_0^{2\sigma}}{M^2}\int_{0}^{\varepsilon_0^{-\frac{\sigma}4}} \rho^{k} d\rho
\leq \frac{ \varepsilon_0^{\frac{2\sigma}3}}{4M^2},$$
and for $k=0,2,4$,
$$
\frac{1}{M^2}\int_{[\inf\sigma(H), \infty)\atop{\rho>\varepsilon_0^{-\frac{\sigma}{4}}}} |f| \left|\partial_{\rho}^2\left( \frac{ \rho^k }{1+\rho^8}\right)\right| \partial\rho \, dE
\leq\frac{|f|_{\Gamma^{(M)}}}{M^2} \int_{\varepsilon_0^{-\frac{\sigma}{4}}}^{\infty} \frac{300\rho^{18}}{(1+\rho^8)^3}  \, d\rho
\leq \frac{\varepsilon_0^{\frac{3\sigma}2}}{4M^2},
$$
\begin{eqnarray*}
\frac{1}{M^2}\int_{[\inf\sigma(H), \infty)\atop{\rho>\varepsilon_0^{-\frac{\sigma}{4}}}}  |\partial h| \left|\partial_{\rho}\left( \frac{ \rho^k }{1+\rho^8}\right)\right|\, dE
&\leq&\frac{\varepsilon_0^{2\sigma}}{M^2} \sum_{l\in\Z^d} c'e^{-\beta'|l|^\iota}
+\frac{2\varepsilon_0^{2\sigma}}{M^2}\int_{\varepsilon_0^{-\frac{\sigma}{4}}}^{\infty} \frac{16\rho^{12}}{(1+\rho^8)^2} d\rho\\
&\leq& \frac{ \varepsilon_0^{\frac{3\sigma}2}}{4M^2}.
\end{eqnarray*}
So (\ref{part0222-1}) is bounded by $\frac{\varepsilon_0^{-\frac{(2k-1)\sigma}8}|f|_{\Gamma^{(M)}}+ \varepsilon_0^{\frac{2\sigma}3}}{2M^2}$, and (\ref{part0222-2}) is bounded by $\frac{\varepsilon_0^{\frac{3\sigma}2}}{2M^2}$.

Combing with (\ref{0ele2-1}) and (\ref{0ele2-2}), we finish the proof of Lemma \ref{lem003}.\qed

By combining Lemma $\ref{lem001}-\ref{lem003}$, for any $f$ which is ${\cal C}^2$ on each connected components of $\Gamma^{(M)}$ and satisfies (\ref{esti_f}), we get that
\begin{eqnarray}
\left|\int_{[\inf\sigma(H), \infty)\atop {\rho\leq \varepsilon_0^{-\frac{\sigma}4}} } f \cdot \rho^k \cos(M \rho)\, \partial \rho\, dE\right|&\leq&\frac{ \varepsilon_{0}^{\frac{\sigma^2}{12}}+ \varepsilon_0^{-\frac{k\sigma}4}|f|_{\Gamma^{(M)}}}{|M|^{1+\frac\sigma{15}}} ,\label{ineq_int_f1}  \\
 \left|\int_{[\inf\sigma(H), \infty)\atop {\rho> \varepsilon_0^{-\frac{\sigma}4}} } \frac{f\cdot \rho^k \cos(M \rho)}{1+\rho^8}\, \partial \rho\, dE\right|&\leq&\frac{ \varepsilon_{0}^{\frac{\sigma^2}{15}}}{|M|^{1+\frac\sigma{15}}}.\label{ineq_int_f2}
\end{eqnarray}
In view of Lemma \ref{coefficient_abc}, we have that, for $l=0,1$, for any $x,\, y\in\R$
$$|\beta_l^{(M)}(x)\beta_l^{(M)}(y)|_{\Gamma^{(M)}}\leq  2\varepsilon_0^{2l\sigma}, \quad |\beta_l^{(M)}(x)\beta_l^{(M)}(y)-\delta_{l,0}|_{\Gamma^{(M)}}\leq  2\varepsilon_0^{2l\sigma}.$$
We apply the inequalities (\ref{ineq_int_f1}) and (\ref{ineq_int_f2}) to $f=h_l=\beta_l^{(M)}(x)\beta_l^{(M)}(y)$ for $k=2l,\, 2l+2$, and to $f=h_l-\delta_{l,0}=\beta_l^{(M)}(x)\beta_l^{(M)}(y)-\delta_{l,0}$ for $k=2l$, since the assumptions of $f$ in (\ref{esti_f}) can be easily deduced from Lemma \ref{coefficient_abc}.
Then Proposition \ref{hmn_integral} is shown.\qed

\

As a application of Proposition \ref{hmn_integral}, we have
\begin{Lemma}\label{coro_hmn_integral}
For $l=0,1$, for any $x,y\in\R$, $M\in\R$ with $|M|>1$, we have
\begin{eqnarray}
\left|\int_{[\inf\sigma(H), \infty)\atop {\rho\leq \varepsilon_0^{-\frac{\sigma}4}} } (\beta_l(x)\beta_l(y)-\delta_{l,0}) \, \rho^{2l} \,  \cos (M \rho) \, \partial \rho\, dE\right| &\leq&\frac{\varepsilon_0^{\frac{\sigma^2}{16}}}{|M|^{1+\frac{\sigma}{15}}},\label{esti_int-1}  \\
\left|\int_{[\inf\sigma(H), \infty)\atop {\rho> \varepsilon_0^{-\frac{\sigma}4}} } \frac{(\beta_l(x)\beta_l(y)-\delta_{l,0}) \,  \rho^{2l}  \cos(M \rho)}{1+\rho^8} \, \partial \rho\, dE\right|&\leq&\frac{\varepsilon_0^{\frac{\sigma^2}{16}}}{|M|^{1+\frac{\sigma}{15}}},\label{esti_int-2}
\end{eqnarray}
and for $k=2l, \, 2l+2$,
\begin{eqnarray}
\left|\int_{[\inf\sigma(H), \infty)\atop {\rho\leq \varepsilon_0^{-\frac{\sigma}4}} } \beta_l(x)\beta_l(y) \, \rho^{k} \,  \cos (M \rho) \, \partial \rho\, dE\right|&\leq&\frac{\varepsilon_0^{\frac{\sigma^2}{16}}+2 \varepsilon_0^{(2l-\frac{k}4)\sigma}}{|M|^{1+\frac{\sigma}{15}}},\label{esti_int-3}  \\
\left|\int_{[\inf\sigma(H), \infty)\atop {\rho> \varepsilon_0^{-\frac{\sigma}4}} } \frac{\beta_l(x)\beta_l(y) \,  \rho^k  \cos(M \rho)}{1+\rho^8} \, \partial \rho\, dE\right|&\leq&\frac{\varepsilon_0^{\frac{\sigma^2}{16}}}{|M|^{1+\frac{\sigma}{15}}}.\label{esti_int-4}
\end{eqnarray}
\end{Lemma}

\proof In (\ref{esti_int-3}) and (\ref{esti_int-4}), we can consider the approximated integrals with respect to $E$ with $\beta_{l}(x) \beta_{l}(y)$ replaced by $\beta^{(M)}_{l}(x) \beta^{(M)}_{l}(y)$, i.e., for $k=2l,\, 2l+2$,
$$\int_{[\inf\sigma(H), \infty)\atop {\rho\leq \varepsilon_0^{-\frac{\sigma}4}}}\beta_l^{(M)}(x)\beta^{(M)}_l(y)\rho^{k} \cos (M\rho)\,\partial \rho \,  dE,$$
$$ \int_{[\inf\sigma(H), \infty)\atop {\rho>\varepsilon_0^{-\frac{\sigma}4}}}\frac{\beta_l^{(M)}(x)\beta^{(M)}_l(y) \rho^{k} \cos (M\rho)}{1+\rho^8}\,\partial \rho \,  dE.$$
Indeed,
$\partial\rho=0$ outside $\sigma(H)$ with $|\sigma(H)\setminus \Sigma|=0$, and, in view of (\ref{error_beta}), $|\beta_l-\beta^{(z)}_l|\leq 10 \varepsilon_J^{\frac14}$ on $\Sigma_j$, $0\leq j\leq J+1$,
with $J=J(M)=\min\{j\in\N:|M|\leq \varepsilon_j^{-\sigma}\}$. Combining with the fact that
$$|\cup_{j\geq J+1} \Sigma_{j+1}|\leq \sum_{j\geq J+1} \ln\varepsilon_0|^{(j+1)^3 d} \varepsilon_{j}^{\sigma}\leq  \varepsilon_{J}^{\sigma+\frac{3\sigma^2}{4}},$$
we can see that the errors are respectively less than
$$20\varepsilon_J^{\frac{1}{4}} \int_0^{\varepsilon_0^{-\frac{\sigma}4}} \rho^{k} \,  d\rho + \frac12\varepsilon_{J}^{\sigma+\frac{\sigma^2}{2}}, \quad  20\varepsilon_J^{\frac{1}{4}} \int_{\varepsilon_0^{-\frac{\sigma}4}}^\infty \frac{ \rho^{k}}{1+\rho^8}  \,  d\rho + \frac12\varepsilon_{J}^{\sigma+\frac{\sigma^2}{2}},$$
both of which are bounded by $\varepsilon_{J}^{\sigma+\frac{\sigma^2}{2}}\leq \frac{\varepsilon_J^{\frac{\sigma}4}}{|M|^{1+\frac\sigma4}}$.

By applying Proposition \ref{hmn_integral} to the the approximated integrals, combining with the errors, we show the estimates for the integrals about $\beta_l(x)\beta_l(y)$. Since it is similar for $\beta_l(x)\beta_l(y)-\delta_{l,0}$ in (\ref{esti_int-1}) and (\ref{esti_int-2}), we finish the proof.
\qed

\section{Modified spectral transformation}\label{section_spec_trans}
\noindent

Let the matrix of measures $d\varphi$ be
\begin{eqnarray*}
  \left . d\varphi \right|_{\Sigma\cap\left\{\rho\leq \varepsilon_0^{-\sigma/4}\right\}} &:=& \frac{1}{\pi}\left(\begin{array}{cc}
                                        (\partial\rho)^{-1} dE & 0 \\[1mm]
                                         0 &  (\partial\rho)^{-1}  dE
                                       \end{array}\right) \\
   \left . d\varphi \right|_{\Sigma\cap\left\{\rho>\varepsilon_0^{-\sigma/4}\right\}} &:=& \frac{1}{\pi}\left(\begin{array}{cc}
                                        \frac{(\partial\rho)^{-1}}{1+\rho^8} dE & 0 \\[1mm]
                                         0 &  \frac{(\partial\rho)^{-1}}{1+\rho^8} dE
                                       \end{array}\right),  \\
  \left. d\varphi\right|_{\R\setminus \Sigma} &:=& 0.
\end{eqnarray*}
Recall the definition of ${\cal L}^2-$space given in (\ref{gen_L2}).
Then ${\cal L}^2(d\varphi)$ means the space of vectors
$G=(g_j)_{j=1,2}$, with $g_j$ functions of $E\in\R$ satisfying
\begin{eqnarray*}
\|G\|_{{\cal L}^2(d\varphi)}^2&:=& \frac{1}{\pi}\int_{\Sigma\cap\left\{\rho\leq\varepsilon_0^{-\sigma/4}\right\}} (|g_1|^2+|g_2|^2)\,  (\partial\rho)^{-1}dE \\
   & & + \frac{1}{\pi}\int_{\Sigma\cap\left\{\rho>\varepsilon_0^{-\sigma/4}\right\}} (|g_1|^2+|g_2|^2)\,  \frac{(\partial\rho)^{-1}}{1+\rho^8}dE<\infty.
\end{eqnarray*}

For the continuous Schr\"odinger operator $H$,
we define the modified spectral transformation ${\cal S}$ on $L^2(\R)$:
$${\cal S} q=\left(\begin{array}{c}
\int_{\R}q(x) {\Cal K}(x) dx \\[1mm]
\int_{\R}q(x) {\Cal J}(x) dx
\end{array}\right), \quad \forall \, q\in L^2(\R),$$
recalling that, on $\Sigma$,
\begin{equation}\label{K_J}
\left(\begin{array}{c}
{\Cal K}(x)\\[1mm]
{\Cal J}(x)
\end{array}\right)=\left(\begin{array}{c}
\beta_0(x) \sin (x\rho)\\[1mm]
\beta_0(x) \cos (x\rho)
\end{array}\right)
+\left(\begin{array}{c}
\beta_1(x) \rho \cos (x\rho)\\[1mm]
- \beta_1(x) \rho \sin (x\rho)
\end{array}\right)=:{\cal B}_0(x)+ {\cal B}_1(x),
\end{equation}
with the coefficients $\beta_l$, $l=0, 1$, satisfying (\ref{esti_abc_sigma_avant}).
%where, for
%$\beta_{l}=\left\{\begin{array}{ll}
%\tilde \beta_{l}, & E\in\Sigma_0 \\[1mm]
%     \xi^{8}\tilde \beta_{l},&  E\in\Sigma_{j+1},\;\ j\geq0
%  \end{array}\right.$,
%with
%$$\tilde \beta_{0}(x) := Y_{11}\left(\omega x\right)B_{12}-Y_{12}\left(\omega x\right)B_{11}, \quad \tilde \beta_{1}(x) :=Y_{12}\left(\omega x\right).$$

\begin{Remark}
The modified spectral transformation ${\cal S}$ is constructed for getting better differentiability with respect to $E$, so it is not necessarily a unitary one.
Comparing with (\ref{classical_uv}) for the free Schr\"odinger operator, ${\cal K}(x)$ and ${\cal J}(x)$ for ${\cal S}$ have a smoothing factor $\xi^{8}$ in a small part of spectrum to cover the singularities.
Moreover, instead of the classical spectral measures shown in Theorem \ref{spectral_measure_matrix},
we use the measures $(\partial\rho)^{-1}dE$ and $\frac{(\partial \rho)^{-1}}{1+\rho^8}dE$ in different parts of spectrum, which has a nice regularity in view of the transversality (\ref{transversality_rho}) of $\partial \rho$ and hence cover the singularity caused by $\partial\rho$.
Moreover, the second measure $\frac{(\partial \rho)^{-1}}{1+\rho^8}dE$ will control the unboundedness of $\rho$.
\end{Remark}

\begin{Remark}
With the purely absolute continuity of the spectrum, we can conclude that the spectral transformation for any non-zero $q\in L^2(\R)$ is supported on a subset of $\sigma(H)$ with positive Lebesgue measure.
Hence, in constructing the modified spectral transformation, we can neglect a zero-measure subset of $\sigma(H)$ and just focus on the full-measure subset $\Sigma$. This is the necessity of the purely absolute continuity in the proof.
\end{Remark}

The following lemma shows that ${\cal S}$ is well defined from $L^2(\R)$ to ${\cal L}^2(d\varphi)$.
\begin{Lemma}\label{well-defined}
Given any $q\in  L^2(\R)$, we have
$\left\|{\cal S}q\right\|_{{\cal L}^2(d\varphi)}< K\| q\|_{L^2(\R)}$ for some $K>0$ independent of $q$.
\end{Lemma}
\proof Let $d\tilde\varphi:=4\rho^2(\partial\rho)^2 d\varphi$, i.e., $\left. d\tilde\varphi\right|_{\R\setminus \Sigma} := 0$,
\begin{eqnarray*}
  \left . d\tilde\varphi \right|_{\Sigma\cap\left\{\rho\leq \varepsilon_0^{-\sigma/4}\right\}} &:=& \frac{4}{\pi}\left(\begin{array}{cc}
                                        \rho^2 \partial\rho \, dE & 0 \\[1mm]
                                         0 & \rho^2 \partial\rho\,  dE
                                       \end{array}\right) \\
   \left . d\tilde\varphi \right|_{\Sigma\cap\left\{\rho>\varepsilon_0^{-\sigma/4}\right\}} &:=& \frac{4}{\pi}\left(\begin{array}{cc}
                                        \frac{\rho^2\partial\rho}{1+\rho^8} dE & 0 \\[1mm]
                                         0 &  \frac{\rho^2\partial\rho}{1+\rho^8} dE
                                       \end{array}\right).
\end{eqnarray*}
To estimate $\|{\cal S}q\|_{{\cal L}^2(d\varphi)}$, we can estimate $\|{\cal S}q\|_{{\cal L}^2(d\tilde\varphi)}$ instead.
Indeed, since $(2\rho)^{-1}< \partial\rho<\infty$ on $\Sigma$, we have
$\|{\cal S}q\|_{{\cal L}^2(d\varphi)}< \|{\cal S}q\|_{{\cal L}^2(d\tilde\varphi)}$.

Given any $q\in L^2(\R)$, we have
$$
\|{\cal S}q\|_{{\cal L}^2(d\tilde\varphi)}\leq\left\|\left(\begin{array}{c}
\int_{\R}q(x){\cal B}^{(1)}_0(x)dx\\[1mm]
\int_{\R}q(x){\cal B}^{(2)}_0(x)dx
\end{array}\right)\right\|_{{\cal L}^2(d\tilde\varphi)}+\left\|\left(\begin{array}{c}
\int_{\R}q(x){\cal B}^{(1)}_1(x)dx\\[1mm]
\int_{\R}q(x){\cal B}^{(2)}_1(x)dx
\end{array}\right)\right\|_{{\cal L}^2(d\tilde\varphi)}.$$
By a direct computation with the formulation in (\ref{K_J}), we have, for $l=0,1$,
\begin{eqnarray}
& &\left\|\left(\begin{array}{c}
\int_{\R}q(x){\cal B}^{(1)}_l(x)dx\\[1mm]
\int_{\R}q(x){\cal B}^{(2)}_l(x)dx
\end{array}\right)\right\|_{{\cal L}^2(d\tilde\varphi)}^2\nonumber\\
&=&  \frac{4}{\pi}\int_{\Sigma\cap\left\{\rho\leq\varepsilon_0^{-\sigma/4}\right\}} \left(\left|\int_{\R}q(x){\cal B}^{(1)}_l(x)dx\right|^2+\left|\int_{\R}q(x){\cal B}^{(2)}_l(x)dx\right|^2\right)  \rho^2\partial\rho\, dE \nonumber\\
 & & + \frac{4}{\pi}\int_{\Sigma\cap\left\{\rho>\varepsilon_0^{-\sigma/4}\right\}} \left(\left|\int_{\R}q(x){\cal B}^{(1)}_l(x)dx\right|^2+\left|\int_{\R}q(x){\cal B}^{(2)}_l(x)dx\right|^2\right)\,  \frac{ \rho^2\partial\rho}{1+\rho^8}dE\nonumber\\
&=&\frac{4}{\pi}
\int_{\R^2} q(x) \bar q(y) \int_{\Sigma\cap\left\{\rho\leq\varepsilon_0^{-\sigma/4}\right\}} \beta_{l}(x) \beta_{l}(y) \rho^{2l+2}
\cos(x-y)\rho \cdot \partial\rho \, dE  \, dx  \, dy\label{integral_L2-1}\\
& &+ \frac{4}{\pi} \int_{\R^2} q(x) \bar q(y) \int_{\Sigma\cap\left\{\rho>\varepsilon_0^{-\sigma/4}\right\}} \frac{ \beta_{l}(x) \beta_{l}(y)  \rho^{2l+2}}{1+\rho^8}
\cos(x-y)\rho \cdot\partial\rho \, dE  \, dx  \, dy,\label{integral_L2-2}
\end{eqnarray}
by noting that $\cos(x\rho)\cos(y\rho)+\sin(x\rho)\sin(y\rho)=\cos(x-y)\rho$.

According to the integration on different regions of $\R^2$, we decompose the integral in (\ref{integral_L2-1}) and (\ref{integral_L2-2}) into ${\cal T}_{l,1}+{\cal T}_{l,2}$ and ${\cal U}_{l,1}+{\cal U}_{l,2}$ respectively, with
\begin{eqnarray*}
{\cal T}_{l,1}&:=&\frac{4}{\pi}\int_{(x,y)\in\R^2\atop{|x-y|\leq\varepsilon_0^{-\frac\sigma8}}}q(x)\bar q(y) \int_{\Sigma\cap\left\{\rho\leq\varepsilon_0^{-\sigma/4}\right\}} \beta_l(x)\beta_l(y)\rho^{2l+2}  \cos(x-y)\rho\, \partial\rho \, dE \, dx \, dy,\\
{\cal T}_{l,2}&:=&\frac{4}{\pi}\int_{(x,y)\in\R^2\atop{|x-y|> \varepsilon_0^{-\frac\sigma8}}}q(x)\bar q(y) \int_{\Sigma\cap\left\{\rho\leq\varepsilon_0^{-\sigma/4}\right\}} \beta_l(x)\beta_l(y)\rho^{2l+2} \cos(x-y)\rho\, \partial\rho \, dE \, dx \, dy,\\
{\cal U}_{l,1}&:=&\frac{4}{\pi}\int_{(x,y)\in\R^2\atop{|x-y|\leq\varepsilon_0^{-\frac\sigma8}}}q(x)\bar q(y) \int_{\Sigma\cap\left\{\rho>\varepsilon_0^{-\sigma/4}\right\}}  \frac{\beta_{l}(x) \beta_{l}(y)  \rho^{2l+2}}{1+\rho^8}  \cos(x-y)\rho \cdot \partial\rho \, dE \, dx \, dy,\\
{\cal U}_{l,2}&:=&\frac{4}{\pi}\int_{(x,y)\in\R^2\atop{|x-y|> \varepsilon_0^{-\frac\sigma8}}}q(x)\bar q(y) \int_{\Sigma\cap\left\{\rho>\varepsilon_0^{-\sigma/4}\right\}}  \frac{\beta_{l}(x) \beta_{l}(y)  \rho^{2l+2}}{1+\rho^8} \cos(x-y)\rho \cdot \partial\rho \, dE \, dx \, dy.
\end{eqnarray*}

By (\ref{esti_abc_sigma_avant}) we have $|\beta_l(x)\beta_l(y)|\leq \frac32$ for any $x,y\in\R$, $l=0,1$. So,
\begin{eqnarray}
|{\cal T}_{l,1}|&\leq& \frac{6\varepsilon_0^{-\frac{5\sigma}4}}{\pi}\int_{(x,y)\in\R^2\atop{|x-y|\leq\varepsilon_0^{-\frac\sigma8}}}|q(x)| |\bar q(y)| dx \, dy \nonumber\\
%&=& \frac{12\varepsilon_0^{-\frac{5\sigma}4}}{\pi}\int_{\R}|q(x)| |\bar q(x+z)| dx  \int_{|z|\leq\varepsilon_0^{-\frac\sigma8}}  dz\\
&\leq&\frac{6\varepsilon_0^{-\frac{5\sigma}4}}{\pi}\int_\R |q(x)|^2 \, dx\int_{|z|\leq\varepsilon_0^{-\frac\sigma8}}  dz\nonumber\\
&\leq&\frac12\varepsilon_0^{-\frac{3\sigma}2}\int_\R |q(x)|^2 \, dx, \label{Tl1}
\end{eqnarray}
\begin{equation}\label{Ul1}
|{\cal U}_{l,1}|\leq \frac{6}{\pi}\int_{(x,y)\in\R^2\atop{|x-y|\leq\varepsilon_0^{-\frac\sigma8}}}|q(x)| |\bar q(y)| dx \, dy \cdot \int_{\varepsilon_0^{-\frac{\sigma}{4}}}^\infty \frac{d\rho}{1+\rho^8}
\leq\varepsilon_0^{\sigma}\int_\R |q(x)|^2 \, dx
\end{equation}
As for ${\cal T}_{l,2}$ and ${\cal U}_{l,2}$, we have
\begin{eqnarray*}
{\cal T}_{l,2}&=&\frac{4}{\pi}\int_{\R}q(x)\bar q(x+z) dx \int_{|z|> \varepsilon_0^{-\frac\sigma8}} \int_{\Sigma\cap\left\{\rho\leq\varepsilon_0^{-\sigma/4}\right\}} \beta_l(x)\beta_l(x+z)\rho^{2l+2} \cos(z\rho) \, \partial\rho \, dE \, dz ,\\
{\cal U}_{l,2}&=&\frac{4}{\pi}\int_{\R}q(x)\bar q(x+z) dx \int_{|z|> \varepsilon_0^{-\frac\sigma8}} \int_{\Sigma\cap\left\{\rho\leq\varepsilon_0^{-\sigma/4}\right\}} \frac{\beta_l(x)\beta_l(x+z)\rho^{2l+2} \cos(z\rho)}{1+\rho^8} \, \partial\rho \, dE \, dz.
\end{eqnarray*}
In view of Lemma \ref{coro_hmn_integral}, there is a constant $K'=K'_{\varepsilon_0}>0$ such that
\begin{eqnarray*}
\left|\int_{\Sigma\cap\left\{\rho\leq\varepsilon_0^{-\sigma/4}\right\}} \beta_l(x)\beta_l(x+z)\rho^{2l+2} \cos(z\rho) \, \partial\rho \, dE\right|&\leq& \frac{K'}{|z|^{1+\frac{\sigma}{15}}},  \\
\left|\int_{\Sigma\cap\left\{\rho\leq\varepsilon_0^{-\sigma/4}\right\}} \beta_l(x)\beta_l(x+z)\rho^{2l+2} \cos(z\rho) \, \partial\rho \, dE\right|&\leq&\frac{K'}{|z|^{1+\frac{\sigma}{15}}}.
\end{eqnarray*}
Hence, for $l=0,1$
$$|{\cal T}_{l,2}|,\, |{\cal U}_{l,2}|\leq\frac{4K'}{\pi}\int_{\R}|q(x)||\bar q(x+z)| dx \int_{|z|> \varepsilon_0^{-\frac\sigma8}} \frac{dz}{|z|^{1+\frac{\sigma}{15}}}\leq \frac{8K'\varepsilon_0^{\frac{\sigma^2}{120}}}{\pi} \int_{\R}|q(x)|^2 dx.$$

Combining with (\ref{Tl1}) and (\ref{Ul1}), we finish the proof.\qed

\smallskip

The following lemma shows that ${\cal S}$ is injective.
\begin{Lemma}\label{injective}
Given any $q\in  L^2(\R)\setminus\{0\}$, we have $\left\|{\cal S}q\right\|_{{\cal L}^2(d\varphi)}>0$.
\end{Lemma}
\proof Let $d\hat\varphi:=(\partial\rho)^2 d\varphi$, i.e., $\left. d\tilde\varphi\right|_{\R\setminus \Sigma} := 0$,
$$ \left . d\hat\varphi \right|_{\Sigma\cap\left\{\rho\leq \varepsilon_0^{-\sigma/4}\right\}} := \frac{1}{\pi}\left(\begin{array}{cc}
                                        \partial\rho \, dE & 0 \\[1mm]
                                         0 & \partial\rho \, dE
                                       \end{array}\right),\;\  \left . d\hat\varphi \right|_{\Sigma\cap\left\{\rho>\varepsilon_0^{-\sigma/4}\right\}} := \frac{1}{\pi}\left(\begin{array}{cc}
                                        \frac{\partial\rho}{1+\rho^8} dE & 0 \\[1mm]
                                         0 &  \frac{\partial\rho}{1+\rho^8} dE
                                       \end{array}\right),$$
%\begin{eqnarray*}
%  \left . d\hat\varphi \right|_{\Sigma\cap\left\{\rho\leq \varepsilon_0^{-\sigma/4}\right\}} &:=& \frac{2}{\pi}\left(\begin{array}{cc}
%                                        \partial\rho \, dE & 0 \\[1mm]
%                                         0 & \partial\rho \, dE
%                                       \end{array}\right), \\
%   \left . d\hat\varphi \right|_{\Sigma\cap\left\{\rho>\varepsilon_0^{-\sigma/4}\right\}} &:=& \frac{2}{\pi}\left(\begin{array}{cc}
%                                        \frac{\partial\rho}{1+\rho^8} dE & 0 \\[1mm]
%                                         0 &  \frac{\partial\rho}{1+\rho^8} dE
%                                       \end{array}\right),  \\
%  \left. d\tilde\varphi\right|_{\R\setminus \Sigma} &:=& 0.
%\end{eqnarray*}
Given any $q\in L^2(\R)\setminus\{0\}$, we have
\begin{equation}\label{tri_ineq}
\left|\|{\cal S}q\|_{{\cal L}^2(d\hat\varphi)}-\left\|\left(\begin{array}{c}
\int_{\R}q(x){\cal B}^{(1)}_0(x)dx\\[1mm]
\int_{\R}q(x){\cal B}^{(2)}_0(x)dx
\end{array}\right)\right\|_{{\cal L}^2(d\hat\varphi)}
\right|\leq \left\|\left(\begin{array}{c}
\int_{\R}q(x){\cal B}^{(1)}_1(x)dx\\[1mm]
\int_{\R}q(x){\cal B}^{(2)}_1(x)dx
\end{array}\right)\right\|_{{\cal L}^2(d\hat\varphi)}.
\end{equation}
By a direct computation, we get, for $l=0,1$,
\begin{eqnarray}
& &\left\|\left(\begin{array}{c}
\int_{\R}q(x){\cal B}^{(1)}_l(x)dx\\[1mm]
\int_{\R}q(x){\cal B}^{(2)}_l(x)dx
\end{array}\right)\right\|_{{\cal L}^2(d\hat\varphi)}^2\nonumber\\
%&=&  \frac{2}{\pi}\int_{\Sigma\cap\left\{\rho\leq\varepsilon_0^{-\sigma/4}\right\}} \left(\left|\int_{\R}q(x){\cal B}^{(1)}_l(x)dx\right|^2+\left|\int_{\R}q(x){\cal B}^{(2)}_l(x)dx\right|^2\right)  \partial\rho\, dE \nonumber\\
% & & + \frac{2}{\pi}\int_{\Sigma\cap\left\{\rho>\varepsilon_0^{-\sigma/4}\right\}} \left(\left|\int_{\R}q(x){\cal B}^{(1)}_l(x)dx\right|^2+\left|\int_{\R}q(x){\cal B}^{(2)}_l(x)dx\right|^2\right)\,  \frac{ \partial\rho}{1+\rho^8}dE\nonumber\\
%&=&\frac{8}{\pi}
%\int_{\R^2} q(x) \bar q(y) \int_{\Sigma\cap\left\{\rho\leq\varepsilon_0^{-\sigma/4}\right\}} \beta_{l}(x) \beta_{l}(y) \cdot \rho^{2l+2}\left( \cos(x\rho)\cos(y\rho)+\sin(x\rho)\sin(y\rho) \right)
%\, \partial\rho \, dE  \, dx  \, dy\\
%& &+ \int_{\R^2} q(x) \bar q(y) \int_{\Sigma\cap\left\{\rho>\varepsilon_0^{-\sigma/4}\right\}} \frac{ \beta_{l}(x) \beta_{l}(y) \cdot \rho^{2l+2}}{1+\rho^8}\left( \cos(x\rho)\cos(y\rho)+\sin(x\rho)\sin(y\rho) \right)
% \, \partial\rho \, dE  \, dx  \, dy\\
&=&\frac{1}{\pi}
\int_{\R^2} q(x) \bar q(y) \int_{\Sigma\cap\left\{\rho\leq\varepsilon_0^{-\sigma/4}\right\}} \beta_{l}(x) \beta_{l}(y) \cdot \rho^{2l}
\cos(x-y)\rho \, \partial\rho \, dE  \, dx  \, dy\label{integral_L2-3}\\
& &+ \frac{1}{\pi} \int_{\R^2} q(x) \bar q(y) \int_{\Sigma\cap\left\{\rho>\varepsilon_0^{-\sigma/4}\right\}} \frac{ \beta_{l}(x) \beta_{l}(y) \cdot \rho^{2l}}{1+\rho^8}
\cos(x-y)\rho \, \partial\rho \, dE  \, dx  \, dy.\label{integral_L2-4}
\end{eqnarray}

For $l=0,1$, we decompose the integrals in (\ref{integral_L2-3}) and (\ref{integral_L2-4}) into ${\cal V}_{l,1}+{\cal V}_{l,2}$ and ${\cal W}_{l,1}+{\cal W}_{l,2}$ respectively, with
\begin{eqnarray*}
{\cal V}_{l,1}&:=&\frac{1}{\pi}\int_{\R^2}q(x)\bar q(y) \int_{\Sigma\cap\left\{\rho\leq\varepsilon_0^{-\sigma/4}\right\}} (\beta_l(x)\beta_l(y)-\delta_{l,0})\rho^{2l}  \cos(x-y)\rho\, \partial\rho \, dE \, dx \, dy \\
&=&\frac{1}{\pi}\int_{\R^2}q(x)\bar q(x+z) \int_{\Sigma\cap\left\{\rho\leq\varepsilon_0^{-\sigma/4}\right\}} (\beta_l(x)\beta_l(x+z)-\delta_{l,0})\rho^{2l}  \cos(z\rho) \, \partial\rho \, dE \, dx \, dz,\\
{\cal V}_{l,2}&:=&\frac{\delta_{l,0}}{\pi}\int_{\R^2}q(x)\bar q(y) \int_{\Sigma\cap\left\{\rho\leq\varepsilon_0^{-\sigma/4}\right\}}  \rho^{2l} \cos(x-y)\rho\, \partial\rho \, dE \, dx \, dy\\
&=&\frac{\delta_{l,0}}{\pi}\int_{\R^2}q(x)\bar q(x+z) \int_0^{\varepsilon_0^{-\frac\sigma4}}  \rho^{2l} \cos(z\rho)  \, d\rho \, dx \, dz,\\
{\cal W}_{l,1}&:=&\frac{1}{\pi}\int_{\R^2}q(x)\bar q(y) \int_{\Sigma\cap\left\{\rho>\varepsilon_0^{-\sigma/4}\right\}}  \frac{(\beta_{l}(x) \beta_{l}(y)-\delta_{l,0}) \rho^{2l}}{1+\rho^8}  \cos(x-y)\rho\, \partial\rho \, dE \, dx \, dy\\
&=&\frac{1}{\pi}\int_{\R^2}q(x)\bar q(x+z) \int_{\Sigma\cap\left\{\rho>\varepsilon_0^{-\sigma/4}\right\}}  \frac{(\beta_{l}(x) \beta_{l}(x+z)-\delta_{l,0}) \rho^{2l}}{1+\rho^8}  \cos(z\rho)\, \partial\rho \, dE \, dx \, dz,\\
{\cal W}_{l,2}&:=&\frac{\delta_{l,0}}{\pi}\int_{\R^2}q(x)\bar q(y) \int_{\Sigma\cap\left\{\rho>\varepsilon_0^{-\sigma/4}\right\}}  \frac{ \rho^{2l}}{1+\rho^8} \cos(x-y)\rho\, \partial\rho \, dE \, dx \, dy\\
&=&\frac{\delta_{l,0}}{\pi}\int_{\R^2}q(x)\bar q(x+z) \int _{\varepsilon_0^{-\frac{\sigma}{4}}}^\infty  \frac{\rho^{2l}}{1+\rho^8} \cos(z\rho) \, d\rho \, dx \, dz.
\end{eqnarray*}

In ${\cal V}_{l,1}$ and ${\cal W}_{l,1}$, $l=0,1$, according to the integration on different regions of $\R^2$, we decompose them respectively as
\begin{eqnarray*}
{\cal V}_{l,1}&=&\frac{1}{\pi}\int_{x\in\R\atop{|z|\leq \varepsilon_0^{-\frac\sigma8}}}q(x)\bar q(x+z) \int_{\Sigma\cap\left\{\rho\leq\varepsilon_0^{-\sigma/4}\right\}} (\beta_l(x)\beta_l(x+z)-\delta_{l,0})\rho^{2l}  \cos(z\rho) \partial\rho \, dE \, dx \, dz\\
& &+\frac{1}{\pi}\int_{x\in\R\atop{|z|> \varepsilon_0^{-\frac\sigma8}}}q(x)\bar q(x+z) \int_{\Sigma\cap\left\{\rho\leq\varepsilon_0^{-\sigma/4}\right\}} (\beta_l(x)\beta_l(x+z)-\delta_{l,0})\rho^{2l}  \cos(z\rho) \partial\rho \, dE \, dx \, dz,\\
{\cal W}_{l,1}&=&\frac{1}{\pi}\int_{x\in\R\atop{|z|\leq \varepsilon_0^{-\frac\sigma8}}}q(x)\bar q(x+z) \int_{\Sigma\cap\left\{\rho>\varepsilon_0^{-\sigma/4}\right\}}  \frac{(\beta_{l}(x) \beta_{l}(x+z)-\delta_{l,0})  \rho^{2l}}{1+\rho^8}  \cos(z\rho) \partial\rho \, dE \, dx \, dz\\
& & +\frac{1}{\pi}\int_{x\in\R\atop{|z|> \varepsilon_0^{-\frac\sigma8}}}q(x)\bar q(x+z) \int_{\Sigma\cap\left\{\rho>\varepsilon_0^{-\sigma/4}\right\}}  \frac{(\beta_{l}(x) \beta_{l}(x+z)-\delta_{l,0}) \rho^{2l}}{1+\rho^8}  \cos(z\rho) \partial\rho \, dE \, dx \, dz.
\end{eqnarray*}
By (\ref{esti_abc_sigma_avant}) we have, for any $x,y\in\R$,
$$|\beta_l(x)\beta_l(y)-\delta_{l,0}|_{\Sigma_0}\leq 2\varepsilon_0^{\frac14}; \quad |\beta_l(x)\beta_l(y)|_{\Sigma_{j+1}}\leq \varepsilon_j^{2\sigma}, \;\ j\geq 0. $$
So, in the region $\{x\in\R, |z|\leq \varepsilon_0^{-\frac\sigma8}\}$, we estimate the integrals with respect to $E$ as
\begin{eqnarray*}
& &\frac{1}{\pi}\left|\int_{\Sigma\cap\left\{\rho\leq\varepsilon_0^{-\sigma/4}\right\}} (\beta_l(x)\beta_l(x+z)-\delta_{l,0})\rho^{2l}  \cos(z\rho) \partial\rho \, dE\right|  \\
&\leq&  \sum_{j\geq 0}\int_{\Sigma_j\cap\left\{\rho\leq\varepsilon_0^{-\sigma/4}\right\}} |\beta_l(x)\beta_l(x+z)-\delta_{l,0}|\rho^{2l}  \partial\rho \, dE \\
&\leq&  2\varepsilon_0^{\frac14}\int_{\Sigma_0\cap\left\{\rho\leq\varepsilon_0^{-\sigma/4}\right\}}\rho^{2l}  \partial\rho \, dE + \sum_{j\geq 0}(1+\varepsilon_j^{2\sigma})\int_{\Sigma_{j+1}\cap\left\{\rho\leq\varepsilon_0^{-\sigma/4}\right\}}\rho^{2l}  \partial\rho \, dE\\
&\leq&  6\varepsilon_0^{\frac14-\frac{3\sigma}4} + 3\sum_{j\geq 0}(1+\varepsilon_j^{2\sigma})\varepsilon_0^{-\frac{\sigma}2} |\ln\varepsilon_0|^{(j+1)^3 d} \varepsilon_{j}^{\sigma}\\
&\leq& \frac12\varepsilon_0^{\frac{3\sigma}8},\\
& &\frac{1}{\pi} \left|\int_{\Sigma\cap\left\{\rho>\varepsilon_0^{-\sigma/4}\right\}}  \frac{(\beta_{l}(x) \beta_{l}(x+z)-\delta_{l,0})  \rho^{2l}}{1+\rho^8}  \cos(z\rho) \partial\rho \, dE \, dx \, dz\right|\\
&\leq&  2\varepsilon_0^{\frac14}\int_{\Sigma_0\cap\left\{\rho>\varepsilon_0^{-\sigma/4}\right\}}\frac{\rho^{2l}  \partial\rho }{1+\rho^8}dE + \sum_{j\geq 0}(1+\varepsilon_j^{2\sigma})\int_{\Sigma_{j+1}\cap\left\{\rho>\varepsilon_0^{-\sigma/4}\right\}}\frac{\rho^{2l}  \partial\rho }{1+\rho^8}dE\\
&\leq& \frac12\varepsilon_0^{\frac{3\sigma}8},
\end{eqnarray*}
recalling the measure of $\Sigma_{j}$ estimated in (\ref{measure_sigma_j}).
In the region $\{x\in\R, |z|> \varepsilon_0^{-\frac\sigma8}\}$, by applying Lemma \ref{coro_hmn_integral}, we get
\begin{eqnarray*}
\left| \int_{[\inf\sigma(H), \infty)\atop {\rho\leq \varepsilon_0^{-\frac{\sigma}4}}}\left(\beta_l(x)\beta_l(x+z)-\delta_{l,0}\right)  \rho^{2l}\cos (z\rho)\,\partial \rho \,  dE\right|&\leq&\frac{\varepsilon_0^{\frac{\sigma^2}{16}}}{|z|^{1+\frac{\sigma}{15}}},  \\
\left|\int_{[\inf\sigma(H), \infty)\atop {\rho>\varepsilon_0^{-\frac{\sigma}4}}}\frac{\left(\beta_l(x)\beta_l(x+z)-\delta_{l,0}\right)  \rho^{2l}\cos (z\rho)}{1+\rho^8}\,\partial \rho \,  dE  \right|&\leq&\frac{\varepsilon_0^{\frac{\sigma^2}{16}}}{|z|^{1+\frac{\sigma}{15}}}.
\end{eqnarray*}
Hence,
\begin{eqnarray}
|{\cal V}_{l,1}|, \, |{\cal W}_{l,1}|&\leq&  \frac{\varepsilon_0^{\frac{\sigma^2}{16}}}{\pi}\int_{\R}|q(x)| |\bar q(x+z)| dx \int_{|z|> \varepsilon_0^{-\frac\sigma8}} \frac{ dz}{|z|^{1+\frac{\sigma}{15}}}\nonumber\\
& & +\frac12\varepsilon_0^{\frac{3\sigma}8} \int_{x\in\R\atop{|z|\leq \varepsilon_0^{-\frac\sigma8}}}|q(x)| |\bar q(x+z)|  dx dz\nonumber \\
&\leq&\frac14\varepsilon_0^{\frac{\sigma^2}{16}} \int_{\R}|q(x)|^2 dx. \label{VW1}
\end{eqnarray}
Obviously, ${\cal V}_{1,2}={\cal W}_{1,2}=0$. So we can deduce that
\begin{equation}\label{qB_1}
\left\|\left(\begin{array}{c}
\int_{\R}q(x){\cal B}^{(1)}_1(x)dx\\[1mm]
\int_{\R}q(x){\cal B}^{(2)}_1(x)dx
\end{array}\right)\right\|_{{\cal L}^2(d\hat\varphi)}={\cal V}_{1,1}+{\cal W}_{1,1}+{\cal V}_{1,2}+{\cal W}_{1,2}\leq \frac12\varepsilon_0^{\frac{\sigma^2}{16}} \int_{\R}|q(x)|^2 dx.
\end{equation}

We still need to consider
\begin{eqnarray*}
{\cal V}_{0,2}&=&\frac{1}{\pi}\int_{\R^2}q(x)\bar q(x+z) \int_0^{\varepsilon_0^{-\frac\sigma4}}  \cos(z\rho)  \, d\rho \, dx \, dz,\\
{\cal W}_{0,2}&=&\frac{1}{\pi}\int_{\R^2}q(x)\bar q(x+z) \int _{\varepsilon_0^{-\frac{\sigma}{4}}}^\infty  \frac{\cos(z\rho)}{1+\rho^8} \, d\rho \, dx \, dz.
\end{eqnarray*}
We have known that $\frac{1}{\pi}\int_{\R}\bar q(x+z) \int_0^{\infty}  \cos(z\rho)  \, d\rho \, dz=\bar q(x)$, then
\begin{equation}\label{VW2}
{\cal V}_{0,2}+{\cal W}_{0,2}=\int_{\R}|q(x)|^2\, dx-\frac{1}{\pi}\int_{\R^2}q(x)\bar q(x+z) \int_{\varepsilon_0^{-\frac\sigma4}}^\infty \frac{\rho^8 \cos(z\rho)}{1+\rho^8}  d\rho \, dx \, dz.
\end{equation}
We are going to show that
\begin{equation}\label{almost-parseval}
\frac{1}{\pi}\left|\int_{\R^2}q(x)\bar q(x+z) \int_{\varepsilon_0^{-\frac{\sigma}{4}}}^\infty \frac{\rho^8 \cos(z\rho)}{1+\rho^8}  d\rho \, dx \, dz\right| \leq \varepsilon_0^{\frac{\sigma}{5}}\int_{\R} |q(x)|^2 dx.
\end{equation}
For the integral $\frac{1}{\pi}\int_{\R}\bar q(x+z) \int_0^{\infty}  \frac{\rho^8 \cos(z\rho)}{1+\rho^8}  d\rho \, dz$, we can see the convergence by Abel's test. So, we can take $R>\varepsilon_0^{-\frac\sigma4}$ large enough such that
$$\frac{1}{\pi}\left|\int_{z> R} [\bar q(x+z)+\bar q(x-z)] \int_{\varepsilon_0^{-\frac\sigma4}}^\infty \frac{\rho^8 \cos(z\rho)}{1+\rho^8}  \, d\rho\, dz \right| \leq \varepsilon_0^{\frac{\sigma}{3}} |q(x)|.$$
Then
\begin{equation}\label{parseval-1}
\frac{1}{\pi}\left|\int_{x\in\R\atop{|z|>R}}q(x)\bar q(x+z) dx \int_{\varepsilon_0^{-\frac{\sigma}{4}}}^\infty \frac{\rho^8 \cos(z\rho)}{1+\rho^8} d\rho \, dz\right|
\leq\varepsilon_0^{\frac{\sigma}{4}}\int_\R |q(x)|^2 dx.
\end{equation}
On the other hand,
\begin{eqnarray}
& & \frac{1}{\pi}\left|\int_{x\in\R\atop{|z|\leq R}}q(x)\bar q(x+z)  \int_{\varepsilon_0^{-\frac{\sigma}{4}}}^\infty  \frac{\rho^8\cos(z\rho)}{1+\rho^8} d\rho\, dx \, dz\right|\nonumber\\
&\leq&\frac{1}{\pi}\int_{\R} |q(x)|^2 dx \cdot \left| \int_{\varepsilon_0^{-\frac{\sigma}{4}}}^\infty \int_{-R}^R \, \cos(z\rho) dz \, \frac{\rho^8}{1+\rho^8} d\rho\right|\nonumber \\
&\leq&\frac{2}{\pi} \int_{\R} |q(x)|^2 dx \left| \int_{\varepsilon_0^{-\frac{\sigma}{4}}}^\infty \frac{\sin(R\rho)}{\rho} \cdot \frac{\rho^8 }{1+\rho^8} d\rho\right|\nonumber\\
&\leq&\frac{2}{\pi}\int_\R |q(x)|^2 dx\left|\left.\frac{\rho^7 \cos(R\rho)}{R(1+\rho^8)}\right|_{\varepsilon_0^{-\frac{\sigma}{4}}}^\infty- \frac{1}{R}\int_{\varepsilon_0^{-\frac{\sigma}{4}}}^\infty \cos(R\rho) \frac{\rho^{14}-7\rho^6 }{(1+\rho^8)^2} d\rho \right|\nonumber\\
&\leq&\frac12\varepsilon_0^{\frac{\sigma}{5}}\int_\R |q(x)|^2 dx.\label{parseval-2}
\end{eqnarray}
So (\ref{almost-parseval}) is shown by combining (\ref{parseval-1}) and (\ref{parseval-2}).
Hence,
\begin{equation}\label{qB_0}
\left|\left\|\left(\begin{array}{c}
\int_{\R}q(x){\cal B}^{(1)}_0(x)dx\\[1mm]
\int_{\R}q(x){\cal B}^{(2)}_0(x)dx
\end{array}\right)\right\|^2_{{\cal L}^2(d\hat\varphi)}- \int_{\R}|q(x)|^2 dx \right|\leq \frac12\varepsilon_0^{\frac{\sigma^2}{16}} \int_{\R}|q(x)|^2 dx.
\end{equation}

In view of (\ref{tri_ineq}), (\ref{qB_1}) and (\ref{qB_0}), for $q\in L^2(\R)\setminus\{0\}$, we have
\begin{equation}\label{almost-unitary}
\left|\|{\cal S}q\|_{{\cal L}^2(d\hat\varphi)}-\|q\|_{L^2(\R)}\right|\leq \varepsilon_0^{\frac{\sigma^2}{16}}\|q\|_{L^2(\R)}.
\end{equation}

Note that the measures $(\partial\rho)^{-1} dE$ and $\frac{(\partial\rho)^{-1}}{{1+\rho^8}} dE$ is absolutely continuous with respect to $\partial\rho\, dE$ and $\frac{\rho^2\partial\rho}{{1+\rho^8}} dE$ respectively, and $\partial \rho$ is positive almost everywhere on $\Sigma$.
We can see that $\|{\cal S}q\|_{{\cal L}^2(d\varphi)}>0$ if $q\in L^2(\R)\setminus\{0\}$.
So we get the conclusion of lemma.\qed

\

For $x\in\R$, ${\cal K}(x)$ and ${\cal J}(x)$ are differentiable in the sense of Whitney on each $\Sigma_j$ and,
by a direct computation, we have
\begin{eqnarray*}
\left(\begin{array}{c}
\partial{\cal B}^{(1)}_0(x)\\[1mm]
\partial{\cal B}^{(2)}_0(x)
\end{array}\right)&=&\left(\begin{array}{c}
(\partial\beta_0)(x) \sin (x\rho)\\[1mm]
(\partial\beta_0)(x) \cos (x\rho)
\end{array}\right)
+\partial\rho\left(\begin{array}{c}
x\beta_0(x) \cos (x\rho)\\[1mm]
-x\beta_0(x)\sin (x\rho)
\end{array}\right),\\
\left(\begin{array}{c}
\partial{\cal B}^{(1)}_1(x)\\[1mm]
\partial{\cal B}^{(2)}_1(x)
\end{array}\right)&=&\rho\, \left(\begin{array}{c}
(\partial\beta_1)(x) \cos (x\rho)\\[1mm]
-(\partial\beta_1)(x) \sin (x\rho)
\end{array}\right)
+ \partial\rho\left(\begin{array}{c}
\beta_1(x) \cos(x\rho)\\[1mm]
-\beta_1(x)\sin(x\rho)
\end{array}\right)\\[1mm]
& & - \, \rho\, \partial\rho\left(\begin{array}{c}
x\beta_1(x) \sin (x\rho)\\[1mm]
x\beta_1(x) \cos (x\rho)
\end{array}\right).
\end{eqnarray*}
Here $\partial  \beta_l(x)$, $l=0,1$, is the derivative in the sense of Whitney on $\Sigma_j$.
Since $\{\Sigma_j\}_{j\geq0}$ are mutually disjoint, $\partial  \beta_l(x)$ and hence $\partial{\cal K}(x)$, $\partial{\cal J}(x)$ are well defined on $\Sigma$.

\begin{Lemma}\label{sobolev-norm}
For any $q\in  L^2(\R)\setminus\{0\}$ with $\|q\|_D< \infty$, we have
\begin{equation}\label{Norm_Diffusion}
\left|\left\| \left(\begin{array}{c}
\int_{\R}q(x) \partial{\Cal K}(x) dx \\[1mm]
\int_{\R}q(x) \partial{\Cal J}(x) dx
\end{array}\right) \right\|_{{\cal L}^2(d\varphi)} - \|q\|_D \right|\leq \varepsilon_0^{\frac{\sigma^2}{20}}\|q\|_D + \varepsilon_0^{\frac{\sigma}4}\|q\|_{L^2(\R)}.
\end{equation}
\end{Lemma}
\proof 
Recall that for $j\geq0$ and $x\in\R$, $|\partial(\beta_l(x)-\delta_{l,0})|_{{\cal C}_W^\nu(\Sigma_j)}\leq \varepsilon_0^{\sigma}$ for $l=0,1$. So
\begin{eqnarray}
& &\left\|\rho^l\left(\begin{array}{c}
\int_{\R}q(x)(\partial\beta_l)(x)\cos(x\rho)dx\\[1mm]
\int_{\R}q(x)(\partial\beta_l)(x)\sin(x\rho)dx
\end{array}\right)\right\|_{{\cal L}^2(d\varphi)}^2\nonumber\\
&\leq&\frac{\varepsilon_0^{2\sigma}}{\pi}
\int_{\R^2}  |q(x)| |\bar q(y)|   \int_{\Sigma\cap\left\{\rho\leq\varepsilon_0^{-\sigma/4}\right\}} \rho^{2l} \, (\partial\rho)^{-1} \, dE \, dx \, dy\nonumber\\
& & + \frac{\varepsilon_0^{2\sigma}}{\pi}
\int_{\R^2}  |q(x)| |\bar q(y)|   \int_{\Sigma\cap\left\{\rho>\varepsilon_0^{-\sigma/4}\right\}} \frac{\rho^{2l}}{1+\rho^8} \, (\partial\rho)^{-1} \, dE \, dx \, dy\nonumber\\
&\leq&\frac{4\varepsilon_0^{2\sigma}}{\pi}
\int_{\R^2}  |q(x)| |\bar q(y)|  \int_0^{\varepsilon_0^{-\sigma/4}} \rho^{2l+2} \,   d\rho \, dx \, dy
 +\frac{4\varepsilon_0^{2\sigma}}{\pi}
\int_{\R^2}  |q(x)| |\bar q(y)|  \int_{\varepsilon_0^{-\sigma/4}}^\infty \frac{\rho^{2l+2}}{1+\rho^8}  \, d\rho \, dx \, dy\nonumber\\
&\leq&\varepsilon_0^{\frac{15\sigma}{16}}
\int_{\R^2}  |q(x)| |\bar q(y)| \, dx \, dy  \nonumber\\
&\leq&\varepsilon_0^{\frac{15\sigma}{16}} \left( \int_{\R}  \frac{ dx}{1+x^2} \right)\left( \int_{\R}  (x^2+1)|q(x)|^2\, dx\right) \nonumber\\
&\leq&\varepsilon_0^{\frac{7\sigma}8}\left(\|q\|_D^2+\|q\|_{L^2(\R)}^2\right), \label{small_1}
\end{eqnarray}
\begin{eqnarray}
& &\left\|\partial\rho\left(\begin{array}{c}
\int_{\R}q(x)\beta_1(x)\cos(x\rho)dx\\[1mm]
\int_{\R}q(x)\beta_1(x)\sin(x\rho)dx
\end{array}\right)\right\|_{{\cal L}^2(d\tilde\varphi)}^2\nonumber\\
&\leq&\frac{\varepsilon_0^{2\sigma}}{\pi}
\int_{\R^2}  |q(x)| |\bar q(y)|   \int_{\Sigma\cap\left\{\rho\leq \varepsilon_0^{-\sigma/4}\right\}} (\partial\rho)\, dE \, dx \, dy\nonumber\\
& & + \frac{\varepsilon_0^{2\sigma}}{\pi}
\int_{\R^2}  |q(x)| |\bar q(y)|   \int_{\Sigma\cap\left\{\rho>\varepsilon_0^{-\sigma/4}\right\}} \frac{\partial\rho}{1+\rho^8}  \, dE \, dx \, dy\nonumber\\
&\leq&\varepsilon_0^{\frac{7\sigma}8}\left(\|q\|_D^2+\|q\|_{L^2(\R)}^2\right).\label{small_2}
\end{eqnarray}
Note that the above bounds still hold if we change the sign of $\beta_l$ or exchange the positions of $\cos(x\rho)$ and $\sin(x\rho)$.

Now, for $\left(\begin{array}{c}
\int_{\R}q(x) \partial{\Cal K}(x) dx \\[1mm]
\int_{\R}q(x) \partial{\Cal J}(x) dx
\end{array}\right)$, it remains the terms
$$\partial\rho\left(\begin{array}{c}
\int_\R xq(x)\beta_0(x) \cos (x\rho) dx \\[1mm]
-\int_\R xq(x)\beta_0(x)\sin (x\rho) dx
\end{array}\right) \quad {\rm and} \quad  \rho\, \partial\rho\left(\begin{array}{c}
\int_\R xq(x)\beta_1(x) \sin (x\rho) dx\\[1mm]
\int_\R xq(x)\beta_1(x) \cos (x\rho) dx
\end{array}\right).$$
It is direct to see that
\begin{eqnarray*}
\left\|\partial\rho\left(\begin{array}{c}
\int_\R xq(x)\beta_0(x) \cos (x\rho) dx \\[1mm]
-\int_\R xq(x)\beta_0(x)\sin (x\rho) dx
\end{array}\right)\right\|_{{\cal L}^2(d\varphi)}&=&\left\|\left(\begin{array}{c}
\int_\R ({\cal X}q)(x){\cal B}^{(1)}_0(x) dx \\[1mm]
-\int_\R ({\cal X}q)(x){\cal B}^{(1)}_0(x) dx
\end{array}\right)\right\|_{{\cal L}^2(d\hat\varphi)},  \\
\left\|\rho\, \partial\rho\left(\begin{array}{c}
\int_\R xq(x)\beta_1(x) \sin (x\rho) dx\\[1mm]
\int_\R xq(x)\beta_1(x) \cos (x\rho) dx
\end{array}\right)\right\|_{{\cal L}^2(d\varphi)}&=&\left\|\left(\begin{array}{c}
\int_\R ({\cal X}q)(x){\cal B}^{(1)}_1(x) dx \\[1mm]
-\int_\R ({\cal X}q)(x){\cal B}^{(1)}_1(x) dx
\end{array}\right)\right\|_{{\cal L}^2(d\hat\varphi)}.
\end{eqnarray*}
where ${\cal X}$ is the multiplication operator $q(x)\mapsto xq(x)$ on the subspace of $\{q\in L^2(\R):\int_\R x^2 |q(x)|^2\, dx<\infty\}$.
In view of (\ref{qB_1}) and (\ref{qB_0}), then we can get
\begin{eqnarray*}
   \left|\left\|\partial\rho\left(\begin{array}{c}
\int_\R xq(x)\beta_0(x) \cos (x\rho) dx \\[1mm]
-\int_\R xq(x)\beta_0(x)\sin (x\rho) dx
\end{array}\right)\right\|_{{\cal L}^2(d\varphi)}- \int_{\R}|q(x)|^2 dx \right|&\leq& \frac12\varepsilon_0^{\frac{\sigma^2}{16}} \int_{\R}|q(x)|^2 dx, \\
\left\|\rho\, \partial\rho\left(\begin{array}{c}
\int_\R xq(x)\beta_1(x) \sin (x\rho) dx\\[1mm]
\int_\R xq(x)\beta_1(x) \cos (x\rho) dx
\end{array}\right)\right\|_{{\cal L}^2(d\varphi)}   &\leq& \frac12\varepsilon_0^{\frac{\sigma^2}{16}} \int_{\R}|q(x)|^2 dx.
\end{eqnarray*}
Together with (\ref{small_1}) and (\ref{small_2}), we can get (\ref{Norm_Diffusion}).\qed

The following lemma shows that $\left(\begin{array}{c}
                   \int_\R q(x)\partial{\cal K}(x) dx \\[1mm]
                   \int_\R q(x)\partial{\cal J}(x) dx
                 \end{array}
\right)$ converges to the derivative of the modified spectral transformation under some suitable conditions.

\begin{Lemma}\label{orderchanging}
For any $q\in L^2(\R)$ with
\begin{itemize}
  \item [(a1)]
  $\left(\begin{array}{c}
                  \int_\R q(x) {\cal K}(x)\,  dx\\[1mm]
                  \int_\R q(x) {\cal J}(x)\,  dx
                 \end{array}
\right)$ convergent to $F=\left(\begin{array}{c}
                   F_1 \\[1mm]
                   F_2
                 \end{array}\right)$ uniformly in $E$,
  \item [(a2)] $\left(\begin{array}{c}
                   \int_\R q(x)\partial{\cal K}(x)\, dx \\[1mm]
                   \int_\R q(x)\partial{\cal J}(x)\, dx
                 \end{array}
\right)$ convergent to $\tilde H=\left(\begin{array}{c}
                   \tilde H_1 \\[1mm]
                   \tilde H_2
                 \end{array}\right)$ in the sense of ${\cal L}^2( d\varphi)$,
\end{itemize}
if $F$ is ${\cal C}^1_W$ on each $\Sigma_j$, then $\partial F= \tilde H$ a.e. on $\Sigma$.
%\begin{itemize}
%  \item [(1)] For any $q\in\ell^2(\Z_+)$, with
%\begin{itemize}
%  \item $\sum_{n\in\Z_+} q_n {\Cal K}_n$ uniformly convergent to $f$,
%  \item $\sum_{n\in\Z_+} q_n \partial{\Cal K}_n$ converges to $\tilde h$ in the sense of $L^2(\R, d\varphi_+)$,
%\end{itemize}
%if $f$ is differentiable in the sense of Whitney on each $\Sigma_j$, then $\partial f= \tilde h$ a.e. .
%  \item [(2)]
%\end{itemize}
\end{Lemma}
\proof For given $x\in\R$ and $l=0,1$, let $\beta_{l}^j$ be the extension of $\beta_{l}$, ${\cal C}^1$ on $[\inf\sigma(H), \infty)$, with $\beta_{l}^j|_{\Sigma_j}=\beta_{l}$, and let
$${\cal K}^j(x):=\beta^j_0(x) \sin (x\rho) - \beta^j_1(x) \rho \cos (x\rho), \quad
  {\cal J}^j(x):=\beta^j_0(x) \cos (x\rho) - \beta^j_1(x) \rho \sin (x\rho).$$
Obviously, ${\cal K}^j(x)$ is absolutely continuous on $[\inf\sigma(H), \infty)$, so for any compactly supported ${\cal C}^1$ function $\phi$ on $[\inf\sigma(H), \infty)$, by the integration by parts,
$$\int_{\Sigma_j} \partial {\Cal K}(x) \cdot \phi\, dE =\left. {\Cal K}(x)  \cdot \phi \right|_{\Sigma_j}  - \int_{\Sigma_j} {\Cal K}(x)  \cdot \partial \phi\, dE.$$
Here $\Sigma_j$ is a Borel set contained in $\sigma(H)$. It can be written as $$\Sigma_j=[\inf\sigma(H), \infty)\setminus \bigcup_{l\geq 0} I_{l},$$ with $\{I_{l}\}_{l\geq 0}$ a sequence of intervals, mutually disjoint, and $\left. {\Cal K}(x)  \cdot \phi \right|_{\Sigma_j}$ is interpreted as
$$\left. {\Cal K}(x)  \cdot \phi \right|_{\Sigma_j}=\left. {\Cal K}^j(x)  \cdot \phi \right|_{[\inf\sigma(H), \infty)}-\sum_{l\geq 0}\left. {\Cal K}^j(x)  \cdot \phi \right|_{I_l}.$$
Since $\beta^j_{l}$, $\phi$ and $\rho$ are all absolutely continuous on $[\inf\sigma(H), \infty)$, we can see the absolute convergence of $\sum_{l\geq 0}\left. {\Cal K}^j(x) \cdot \phi \right|_{I_l}$. Hence, by Fubini's theorem,
$$\int_{\R} q(x)  (\left.{\Cal K}(x) \cdot  \phi \right|_{\Sigma_j}) dx=\left.\left(\int_{\R} q(x) {\Cal K}(x) dx \right)\cdot  \phi \right|_{\Sigma_j}= \left. F_1 \cdot \phi \right|_{\Sigma_j}.$$
On the other hand, for each $\Sigma_j$, we have, by (a2),
\begin{eqnarray*}
& & \int_{\Sigma_j} |\int_{|x|\leq N} q(x) \partial{\Cal K}(x) \, dx  - \tilde H_1 |\cdot |\phi| dE\\
&\leq&\left(\int_{\Sigma_j\atop{\rho>\varepsilon_0^{-\sigma/4}}} |\int_{|x|\leq N} q(x) \partial{\Cal K}(x)dx - \tilde H_1|^2 \frac{(\partial\rho)^{-1}dE}{1+\rho^8}\right)^{\frac12} \left(\int_{\Sigma_j\atop{\rho>\varepsilon_0^{-\sigma/4}}} |\phi|^2(1+\rho^8)\partial\rho dE\right)^{\frac12}\\
& &+ \, \left(\int_{\Sigma_j\atop{\rho\leq\varepsilon_0^{-\sigma/4}}} |\int_{|x|\leq N} q(x) \partial{\Cal K}(x)dx - \tilde H_1|^2 (\partial\rho)^{-1}\, dE\right)^{\frac12} \left(\int_{\Sigma_j\atop{\rho\leq\varepsilon_0^{-\sigma/4}}} |\phi|^2\partial\rho\, dE\right)^{\frac12} ,
\end{eqnarray*}
which goes to $0$ as $N\rightarrow\infty$.
Hence,
\begin{eqnarray*}
\int_{\Sigma_j} \partial F_1\cdot\phi\, dE &=&\left. F_1 \cdot \phi \right|_{\Sigma_j}  - \int_{\Sigma_j} F_1 \cdot \partial \phi\, dE \\
&=&\int_\R q(x) (\left.{\Cal K}(x)\cdot  \phi \right|_{\Sigma_j}) dx-  \int_\R q(x)\int_{\Sigma_j} {\Cal K}(x) \cdot \partial \phi\, dE \, dx,\\
&=&\int_\R q(x) \int_{\Sigma_j} \partial {\Cal K}(x) \cdot \phi\, dE\, dx\\
&=&\lim_{N\rightarrow \infty} \int_{\Sigma_j} \int_{|x|\leq N} q(x) \partial {\Cal K}(x) \cdot \phi\, dE \, dx \\
&=&\int_{\Sigma_j} \tilde H_1 \cdot \phi\, dE.
\end{eqnarray*}
%\begin{eqnarray*}
%\int_{\Sigma_j} \partial f\cdot\phi dE&=& \left. f \cdot \phi \right|_{\Sigma_j}  - \int_{\Sigma_j} f \cdot \partial \phi\, dE\\
%&=& \left. \sum_n q_n {\Cal K}_n\cdot  \phi \right|_{\Sigma_j}-  \int_{\Sigma_j} \sum_n q_n {\Cal K}_n \cdot \partial \phi\, dE\\
%&=&  \sum_n q_n \int_{\Sigma_j} \partial {\Cal K}_n \cdot \phi\, dE\\
%&=& \lim_{N\rightarrow \infty} \int_{\Sigma_j} \sum_{n\leq N} q_n \partial {\Cal K}_n \cdot \phi\, dE\\
%&=& \int_{\Sigma_j} \tilde h \cdot \phi\, dE.
%\end{eqnarray*}
So $\partial F_1=\tilde H_1$ a.e. on each $\Sigma_j$, hence a.e. on $\Sigma$. Similarly, $\partial F_2=\tilde H_2$ a.e. on $\Sigma$. \qed

\section{Linear growth of the diffusion norm}\label{proof_theorem}
\noindent

Now, let $q(t)$ be the solution to the dynamical equation ${\rm i}\dot{q}=H q$, with $q(0)\in H^{1}(\R)$. Let $G(E,t):=({\cal S}q)(E,t)$. Since, for any $E\in\Sigma$,
$$\frac{1}{\delta}(G(E,t+\delta)-G(E,t))=\frac{1}{\delta}\left(\begin{array}{c}
\int_{\R} [q(x,t+\delta)-q(x,t)] {\cal K}(x) dx \\[1mm]
\int_{\R} [q(x,t+\delta)-q(x,t)] {\cal J}(x) dx
\end{array}\right) \ {\rm for} \ \delta >0,$$
we can verify the differentiability of $G(E,t)$ with respect to $t$.
Then, for $E\in\Sigma$,
\begin{eqnarray*}
{\rm i}\partial_t G(E,t)&=&\left(\begin{array}{c}
\int_{\R} (Hq)(x,t) {\cal K}(E,x) dx \\[1mm]
\int_{\R} (Hq)(x,t) {\cal J}(E,x) dx
\end{array}\right) \\[1mm]
&=&\left(\begin{array}{c}
\int_{\R}  q(x,t) [-\partial^2_x{\cal K}(E,x)+V(\omega x){\cal K}(E,x) ]dx\\[1mm]
\int_{\R}  q(x,t) [-\partial^2_x{\cal J}(E,x)+V(\omega x){\cal J}(E,x) ]dx
\end{array}\right)  \\[1mm]
&=&E G(E,t) .
\end{eqnarray*}
so $G(E,t)= e^{-{\rm i}Et} G(E,0)$.
%Hence, if $G(E,0)$ is differentiable with respect to $E$, then
%\begin{equation}\label{derivative_g}
%\partial G(E,t)=-{\rm i} t\cdot e^{-{\rm i}Et} G(E,0) + e^{-{\rm i}Et}\partial G(E,0).
%\end{equation}

\begin{Corollary}\label{sum_derivation_exchange}
For any solution $q(t)$ to the equation ${\rm i}\dot{q}=H q$, with $q(0)$ compactly supported, we have, for a.e. $E\in\Sigma$,
  \begin{equation}\label{verif_changing_whole}
  \left(\begin{array}{c}
                   \int_{\R} q(x,t)\partial{\cal K}(E,x) dx \\[1mm]
                    \int_{\R} q(x,t)\partial{\cal J}(E,x) dx
                 \end{array}
\right) = -{\rm i} t\cdot e^{-{\rm i}Et} G(E,0) + e^{-{\rm i}Et}\partial G(E,0).
  \end{equation}
\end{Corollary}
\proof Since $q(0)$ is supported on a compact subset of $\R$, saying $\Lambda\subset\R$, we have that $\partial G(E,0)$ is well defined on each $\Sigma_j$, with
$$\partial G(E,0)=\left(\begin{array}{c}
                   \int_{\Lambda} q(x, 0)\partial{\cal K}(x,E) dx \\[2mm]
                    \int_{\Lambda} q(x,0)\partial{\cal J}(x,E) dx
                 \end{array}
\right).$$
Hence, $G(E,t)= e^{-{\rm i}Et} G(E,0)$ is differentiable in the sense of Whitney on each $\Sigma_j$, with
$$\partial G(E,t)=-{\rm i} t\cdot e^{-{\rm i}Et} G(E,0) + e^{-{\rm i}Et}\partial G(E,0).$$
For any finite $t$, $\int_{\R} x^2 |q(x,t)|^2dx<\infty$, which implies $\int_{\R}|q(x,t)|dx<\infty$.
Then, recalling the expression of ${\cal K}(x)$ and ${\cal J}(x)$, we have
$$\int_{\R} |q(x,t) {\cal K}(x)| dx, \; \int_{\R} |q(x,t) {\cal J}(x)| dx \leq (2+\varepsilon_0^{\sigma}\rho)\int_{\R} |q(x,t)| dx, $$
and by Lemma \ref{sobolev-norm}, for $N>0$ sufficiently large,
\begin{eqnarray*}
\left\|\left(\begin{array}{c}
                   \int_{|x|>N} q(x,t)\partial{\cal K}(x) dx\\[2mm]
                    \int_{|x|>N} q(x,t)\partial{\cal J}(x) dx
                 \end{array}
\right)\right\|_{{\cal L}^2(d\varphi)}&\leq& \frac32\left(\int_{|x|>N}x^2|q(x,t)|^2 dx\right)^{\frac12} \\
& &+\, \varepsilon_0^{\frac\sigma4}\left(\int_{|x|>N} |q(x,t)|^2 dx\right)^{\frac12}.
\end{eqnarray*}
So the assumptions (a1) and (a2) of Lemma \ref{orderchanging} are verified.
Applying Lemma \ref{orderchanging}, the proof of (\ref{verif_changing_whole}) is finished.\qed

Given any solution $q(t)$ to ${\rm i}\dot{q}=H q$ with initial datum $q(0)\in H^{1}(\R)$ with $\|q(0)\|_D< \infty$.
Let $\{q^N(0)\}_{N\in\Z_+}\subset H^1(\R)$ be the sequence of functions satisfying that
$$\lim_{N\rightarrow\infty}\|q^N(0)-q(0)\|_{H^1(\R)}= \lim_{N\rightarrow\infty}\|q^N(0)-q(0)\|_{D}=0,$$
with each $q^N(0)$ compactly supported on $[-N,N]$ satisfying $\|q^N(0)\|_{L^2(\R)}\leq\|q(0)\|_{L^2(\R)}$.
By Lemma \ref{well-defined},
$$\lim_{N\to\infty}\|{\cal S} q^N(0)\|_{{\cal L}^2(d\varphi)}=\|{\cal S} q(0)\|_{{\cal L}^2(d\varphi)}.$$
Let $q^N(t)$ be the solution satisfying ${\rm i}\dot{q}^N=H q^N$ with initial datum $q^N(0)$,
and $G_N(E,t)=({\cal S}q^N)(E,t)$.

In view of Lemma \ref{sobolev-norm} and Corollary \ref{sum_derivation_exchange}, we can see,
$$
\left|t\|G_N(E,0) \|_{{\cal L}^2(d\varphi)}- \|q^N(t)\|_D\right|
\leq \varepsilon_0^{\frac{\sigma^2}{20}} \|q^N(t)\|_D+ \varepsilon_0^{\frac{\sigma}4}\|q^N(t)\|_{L^2(\R)} + \|\partial G_N(E,0) \|_{{\cal L}^2(d\varphi)}.
$$
%On the other hand, by (\ref{derivative_g}),
%\begin{equation}\label{truncation_property_2}
%\left|\|\partial G_N(E,t) \|_{{\cal L}^2(d\varphi)}-t\|G_N(E,0)\|_{{\cal L}^2(d\varphi)}\right|\leq \|\partial G_N(E,0) \|_{{\cal L}^2(d\varphi)}.
%\end{equation}
%Combining (\ref{truncation_property_1}) and (\ref{truncation_property_2}),
Hence, we have
$$\frac{\|G_N(E,0)\|_{{\cal L}^2(d\varphi)}-t^{-1}{\cal G}_N(t)}{1+\varepsilon_0^{\frac{\sigma^2}{20}}}\leq t^{-1}\|q^N(t)\|_D\leq \frac{\|G_N(E,0)\|_{{\cal L}^2(d\varphi)}+t^{-1}{\cal G}_N(t)}{1-\varepsilon_0^{\frac{\sigma^2}{20}}}$$
with ${\cal G}_N(t):= \varepsilon_0^{\frac{\sigma}4}\|q^N(t)\|_{L^2(\R)}+\|\partial G_N(E,0) \|_{{\cal L}^2(d\varphi)}$.
By Lemma \ref{sobolev-norm} and Corollary \ref{orderchanging}, combining with the $L^2-$conservation, we can see
$${\cal G}_N(t)\leq 2\varepsilon_0^{\frac{\sigma}4}\|q^N(0)\|_{L^2(\R)} + 2\|q^N(0)\|_D\leq 2\varepsilon_0^{\frac{\sigma}4} \|q(0)\|_{L^2(\R)}+2\|q(0)\|_{D}.$$
So, for $t$ large enough(independent of $N$), $t^{-1}{\cal G}_N(t)$ goes to zero, and
\begin{equation}\label{growth_finite-supp_whole}
\frac{\|{\cal S} q^N(0)\|_{{\cal L}^2(d\varphi)}}{1+\varepsilon_0^{\frac{\sigma^2}{24}}}\leq t^{-1}\|q^N(t)\|_D \leq \frac{\|{\cal S}q^N(0)\|_{{\cal L}^2(d\varphi)}}{1-\varepsilon_0^{\frac{\sigma^2}{24}}}.
\end{equation}
By the ballistic upper bound (\ref{upper_bound}), we have
$$\lim_{N\to\infty}t^{-1}\|q^N(t)-q(t)\|_D\leq c\lim_{N\to\infty} \left(\|q^N(0)-q(0)\|_{D}+\|q^N(0)-q(0)\|_{H^1(\R)}\right)=0.$$
Combining with the fact that $\lim_{N\to\infty}\|{\cal S} q^N(0)-{\cal S} q(0)\|_{{\cal L}^2(d\varphi)}=0$, we can pass
(\ref{growth_finite-supp_whole}) to $N\to\infty$.
Then Theorem \ref{thm-qp} is proven with
$$C=\|{\cal S} q(0)\|_{{\cal L}^2(d\varphi)}, \quad \zeta=\frac{\sigma^2}{24}.$$

\appendix

\section{Fourier transform and ballistic motion for the Laplacian}\label{fourier}
\noindent

Given $q=q(\xi)\in L^1(\R)$, the Fourier transform is given by
$$q(\xi)\mapsto \hat q(x):=\int_\R e^{-2\pi{\rm i}\xi x} q(\xi)\, d\xi.$$
If $q\in H^1(\R)$, then $\hat { q' }(x)=2\pi {\rm i} x \, \hat q(x)$.

We consider the time dependent equation ${\rm i}\partial_t \hat q(x,t)=-\partial^2_x\hat  q(x,t)$. By the Fourier transform, it is transformed into the linear equation
${\rm i}\partial_t q(\xi,t)= 4\pi^2 \xi^2 q(\xi,t)$.
We can easily get the solution
$
q(\xi,t)=e^{- 4{\rm i}\pi^2 \xi^2 t} q(\xi,0)
$,
which implies
\begin{equation}\label{Fourier_growth}
\partial_\xi q(\xi,t)=-8{\rm i}\pi^2 \xi t e^{- 4{\rm i}\pi^2 \xi^2 t} q(\xi,0)+e^{- 4{\rm i}\pi^2 \xi^2 t} \partial_\xi q(\xi,0).
\end{equation}
On the other hand, by Parseval's identity, we have
\begin{equation}\label{Fourier_norm}
\int_\R  |\partial_\xi q(\xi,t)|^2 d\xi=4\pi^2\int_\R x^2 |\hat q(x,t)|^2 dx.
\end{equation}
We can conclude the linear growth of $\|\hat q(t)\|_D$ by combining (\ref{Fourier_growth}) and (\ref{Fourier_norm}).

\

\noindent {\bf Acknowledgements.}
The author would like to thank H. Eliasson and C. Chavaudret for many fruitful discussions about this work. He would also thank L. Stolovitch for supporting this work.
The author would also thank the anonymous referee for valuable comments and suggestions on the manuscript.

\end{document}